\documentclass[
    12pt,
    reqno,
    final,
    ]{amsart}

\usepackage{amsmath, 
    amsthm, 
    amssymb, 
    amsthm, 
    wasysym, 
    enumitem,
    mathrsfs,
}
\usepackage{graphicx}
\usepackage[top             =2cm, 
            bottom          =4.5cm, 
            left            =2.5cm, 
            right           =2.5cm, 
            marginparwidth  =2.0cm,
            ]{geometry}
\usepackage{xcolor}
\usepackage[urlcolor    = blue,
            ]{hyperref}
\usepackage[nameinlink]{cleveref}
\usepackage{enumitem}
\usepackage{tikz}
\usetikzlibrary{patterns}
\usetikzlibrary{cd}
\usepackage{soul}
\usepackage{bbm}
\usepackage{pgffor}
\usepackage{thmtools}
\usepackage{stackengine}
\usepackage{import}
\usepackage[
    colorinlistoftodos,
    prependcaption,
    textsize=tiny,
    color=cyan!50,
    obeyDraft,
    ]{todonotes}
\setuptodonotes{inline}
\usepackage{xargs}
\usepackage{mathtools}
\usepackage{ifdraft}
\usepackage{float}

\usepackage{color}   

\definecolor{darkish-green}{HTML}{32CD32}
\hypersetup{
    colorlinks=true, 
    linktoc=all,     
    linkcolor=blue,  
    citecolor=darkish-green,
}

\setlist[enumerate]{label=(\arabic*)}
\theoremstyle{plain}
\newtheorem{theorem}{Theorem}[section]
\newtheorem{lemma}[theorem]{Lemma}
\newtheorem{proposition}[theorem]{Proposition}
\newtheorem{corollary}[theorem]{Corollary}
\newtheorem{conjecture}[theorem]{Conjecture}

\theoremstyle{definition}
\newtheorem{definition}[theorem]{\textit{Definition}}

\newtheorem{example}[theorem]{\textit{Example}}

\AtBeginEnvironment{example}{%
  \pushQED{\qed}%
}
\AtEndEnvironment{example}{\popQED\endexample}

\theoremstyle{remark}

\newtheorem*{note}{\bf{Note}}

\numberwithin{equation}{section}
\numberwithin{claim}{theorem}

\newenvironment{sketch}{%
  \proof
  
  }{\endproof}

\newcommand{\term}[1]{\textit{#1}}

\newcommand{\set}[1]{\left\{#1\right\}}

\DeclarePairedDelimiter\gen\langle\rangle
\newcommand{\sothat}{\;\middle|\;}

\newcommand{\sse}{\subseteq}

\newcommand{\GG}{\mathbb{G}}
\newcommand{\HH}{\mathbb{H}}

\newcommand{\NN}{\mathbb{N}}

\newcommand{\QQ}{\mathbb{Q}}
\newcommand{\RR}{\mathbb{R}}

\newcommand{\TT}{\mathbb{T}}

\newcommand{\ZZ}{\mathbb{Z}}

\newcommand{\Heis}{\mathcal{H}}

\newcommand{\im}[1]{\operatorname{im}\l(#1\r)}

\newcommand{\biject}{\leftarrow \hspace{-.5em} \rightarrow}
\newcommand{\onetoone}{\biject}
\newcommand{\into}{\hookrightarrow}

\newcommand{\mf}[1]{\mathfrak{#1}}
\newcommand{\mc}[1]{\mathcal{#1}}
\newcommand{\ms}[1]{\mathscr{#1}}

\newcommand{\Id}{\mathbbm{1}}

\renewcommand{\tilde}{\widetilde}

\renewcommand{\bar}{\overline}

\renewcommand{\l}{\left}
\renewcommand{\r}{\right}

\usepackage[OT2,T1]{fontenc}
\DeclareSymbolFont{cyrletters}{OT2}{wncyr}{m}{n}
\DeclareMathSymbol{\Sha}{\mathalpha}{cyrletters}{"58}

\newcommand{\dirlim}{\varinjlim}
\newcommand{\invlim}{\varprojlim}

\newcommand{\dist}{\mathrm{dist}}
\newcommand{\isom}{\cong}

\newcommand{\Hull}[1]{\Omega_{#1}}

\newcommand{\boundary}{\partial}

\renewcommand{\phi}{\varphi}

\newenvironment{diagram}{\begin{center}\begin{tikzcd}[ampersand replacement=\&]}{\end{tikzcd}\end{center}}

\makeatletter
\providecommand\@dotsep{5}
\renewcommand{\listoftodos}[1][\@todonotes@todolistname]{%
  \@starttoc{tdo}{#1}}
\makeatother

\newcommandx{\toprove}[2][1=]{\todo[inline,linecolor=red,backgroundcolor=red!50,bordercolor=black,#1]{Missing proof. #2}}
\newcommandx{\todefine}[2][1=]{\todo[inline,linecolor=red,backgroundcolor=red!50,bordercolor=black,#1]{Missing definition. #2}}

\newcommandx{\tofinish}[2][1=]{\todo[inline,linecolor=yellow,backgroundcolor=yellow!50,bordercolor=black,#1]{Incomplete. #2}}

\newcommandx{\tofix}[2][1=]{\todo[inline,linecolor=orange,backgroundcolor=orange!50,bordercolor=black,#1]{Errors present. #2}}

\newcommandx{\toreview}[2][1=]{\todo[inline,linecolor=green,backgroundcolor=green!50,bordercolor=black,#1]{Review. #2}}

\newcommandx{\toedit}[2][1=]{\todo[inline,linecolor=blue,backgroundcolor=blue!25,bordercolor=black,#1]{Make edits #2}}

\newcommandx{\tocite}[2][1=]{\todo[inline,linecolor=black,backgroundcolor=gray!20,bordercolor=black,#1]{Citation needed #2}}

\newcommandx{\question}[2][1=]{\todo[inline,linecolor=black,backgroundcolor=white,bordercolor=red,textcolor=red,#1]{\underline{Question:} #2}}

\newcommandx{\answer}[2][1=]{\todo[inline,linecolor=black,backgroundcolor=white,bordercolor=green,textcolor=green,#1]{\underline{Answer:} #2}}

\newcommandx{\starthere}[2][1=]{\todo[inline,linecolor=red,backgroundcolor=red!50,bordercolor=black,#1]{\textbf{START HERE!}}}

\newcommand{\rsimplex}[1][\Delta]{{#1}}
\newcommand{\eusimplex}{\boldsymbol{\Delta}}
\newcommand{\Del}{\operatorname{Del}}
\newcommand{\Conv}{\operatorname{Conv}}

\newcommand{\tube}[2]{#1 #2}

\newcommand{\pathspace}[1]{P{#1}^{(1)}}
\newcommand{\loopspace}[1]{\Omega{#1}^{(1)}}
\graphicspath{{../07 Figures/}}

\title{Simple Tilings of Riemannian Homogeneous Spaces}
\author{Kyle Hansen}

\begin{document}
\markboth{Hansen}{\thetitle}

\begin{abstract}
    We define simple tilings in the general context of a $G$-tiling on a Riemannian homogeneous space $\mc{M}$ to be tilings by Riemannian simplices.
As evidence that this definition is natural, we prove that a large class of tilings of $\mc{M}$ are MLD to simple ones.
We further show that there is an isomorphism between \u{C}ech cohomology and pattern-equivariant cohomology of simple tilings of homogeneous Riemannian Spaces.
This demonstrates the utility of our definition by generalizing an established result on the topology of tiling spaces of \texorpdfstring{$\RR^n$}{Euclidean space}.
\end{abstract}

\maketitle


\section{Introduction}\label{sec:intro}
    Any sufficiently nice tiling $\mc{T}$ of $\RR^n$ can be modeled using a \textit{simple tiling}, that is, a tiling by convex polytopes meeting full-face to full-face.
To create such a simple tiling (see \cite[\S 4.2]{priebe2000towards}), start by identifying an indicator point of each tile, for instance a center of mass.
The resulting Delaunay set $D \subseteq \RR^n$ may be studied in its own right as a model of $\mc{T}$, though without the additional geometric data of a particular tiling.
To create the desired simple tiling, one constructs a Vorono\"{i} diagram $\mc{T}'$ from $D$.
The patterns of $\mc{T}'$ are indistinguishable from those of $\mc{T}$, so that the two tilings are \textit{mutually locally derivable} (MLD), the strongest notion of equivalence between tilings.
The study of $\mc{T}$ then reduces to the study of the simple tiling $\mc{T}'$, which by design contains elegant geometric and combinatorial structure.

Unfortunately, the definition of a simple tiling in $\RR^n$ does not extend immediately to more general contexts: it is unreasonable, for instance, to assume that tiles are convex polytopes in non-Euclidean spaces such as the (real) Heisenberg group.
While the study of the intermediate Delaunay set $D$ has progressed beyond the Euclidean setting (see \cite{machado2018approximate}), to the author's knowledge there has not yet been a formalization of simple tilings developed for non-Euclidean spaces.
In this paper, we generalize the definition of a simple tiling to be one whose tiles are \term{Riemannian simplices}. 
Properties of such simplices are developed in \cite{dyer2015riemannian}, and we define them precisely in \autoref{sec:delaunay}. 
For now it suffices to say that in spaces of bounded geometry, at a small enough scale, such simplices are determined by their vertex set and the local geometry of the underlying space. 
This makes them a natural analogue of the standard Euclidean simplex, hence of convex polyhedra, in non-Euclidean settings. 
To demonstrate the naturality of this definition, we prove the following:

\begin{restatable*}{theorem}{simpleTilings}\label{thm:non-Euclidean normal tilings are simple}
    Let $\mc{M}$ be a Riemannian homogeneous space and $G \le \mathrm{Isom}(\mc{M})$ acting transitively on $\mc{M}$.
    Given a geometrically normal $G$-tiling $\mc{T}$ of $\mc{M}$, there is a simple tiling $\mc{T}_{\Delta}$ which is MLD to $\mc{T}$. 
\end{restatable*}

A \textit{geometrically normal} tiling is one whose tiles have geometrically tame boundaries, which only intersect in reasonable ways. 
We define this notion precisely in \autoref{sec:tiling spaces}.

\autoref{thm:non-Euclidean normal tilings are simple} can be applied, for example, to $G = \mc{M}$ a connected, simply connected, nilpotent Lie group acting on itself by left-translation.
Such spaces naturally generalize $\RR^n$, and the theory of mathematical quasicrystals naturally extends to this setting as well.
In recent years (see \cite{bjorklund2018approximate,machado2018approximate,kaiser2022complexity} for example), there has been a growing interest in the study of mathematical quasicrystals in these spaces.

The general outline of the proof of \autoref{thm:non-Euclidean normal tilings are simple} is rather straightforward. 
Begin by putting a sufficiently nice indicator set of points on a tubular neighborhood of each prototile of $\mc{T}$. 
Do this in a way that is both consistent between admissible adjacencies of prototiles (i.e. when gluing prototiles along their faces, the point sets on the overlap of their neighborhoods align), and in such a way that their union is a Delaunay set $D$. 
Ensuring that $D$ is generic will induce a triangulation of $\mc{M}$ by Riemannian simplices, whose vertex set is $D$ itself.

The main technical difficulty one must overcome is finding a way to ensure that $D$ is in fact a generic Delaunay set. 
The first approach one might take is to simply stick a generic Delaunay set on each prototile; however, because the Delaunay set must be generic when tiles are laid next to one another, it is possible that this approach creates degeneracies (i.e., places where the Delaunay set is \textit{not} generic) across boundaries of tiles.

Moreover, one must take into account \textit{all} the possible adjacencies which might simultaneously affect the genericity of a Delaunay set on a particular tile.
We do so by considering the \textit{Anderson-Putnam complex} $\Gamma$ of $\mc{T}$, which encodes all relevant adjacencies of each tile.
We proceed inductively up tubular neighborhoods of the skeleta of $\Gamma$.
The algorithms and analysis of the series of papers \cite{boissonnat2013stability,boissonnat2014delaunay,boissonnat2016obstruction,boissonnat2018delaunay} provide the tools necessary to create a Delaunay set on $\Gamma$ which pulls back to a generic, pattern-equivariant Delaunay set on $\mc{T}$.

Having established the naturality of this definition of a simple tiling, we expect that many theorems of simple tilings for $\RR^n$ will extend to simple tilings of more general spaces.
We demonstrate one such result:

\begin{restatable*}
    {theorem}{PECohomology}\label{thm:PE cohomology}
    Given a simple $G$-tiling $\mc{T}$ of a Riemannian homogeneous space $\mathcal{M}$, the (real or integer) pattern-equivariant cohomology of $\mc{T}$ is isomorphic to the \u{C}ech cohomology of the hull $\Hull{\mc{T}}$ of $\mc{T}$ (with real or integer coefficients, respectively).
\end{restatable*}

In \autoref{sec:tiling spaces}, we present the relevant tools for studying tiling spaces and review some of the well-known results for simple tilings of $\RR^n$, as well as examples of tilings more generally.
In \autoref{sec:delaunay}, we discuss Delaunay sets in detail and state a few necessary results about Delaunay triangulations of manifolds.
In \autoref{sec:simple tilings}, we generalize the definition of simple tilings of homogeneous Riemannian spaces, and apply the previous two sections to prove \autoref{thm:non-Euclidean normal tilings are simple}. 
In \autoref{sec:pattern equivariant}, we prove \autoref{thm:PE cohomology}, demonstrating that our definition of "simple tilings" is a useful generalization.

We furthermore include an appendix investigating the topological structure of tiling spaces of rational, connected, simply connected, nilpotent Lie groups.
In \autoref{sec:nilpotent lie groups}, we provide the necessary background for such groups.
In \cite[Theorem 1]{sadun2003tilingfiberbundles} it is demonstrated that tiling spaces of $\RR^n$ are fiber bundles over a torus; in \autoref{sec:fiber bundles}, we demonstrate progress towards extending this result to tiling spaces of rational, connected, simply connected, nilpotent Lie groups.
Further directions of ongoing and future work are presented in \autoref{sec:further directions}.

\subsubsection{Acknowledgements} 
The author would like to thank his advisor, Fedya Manin, for inspiration and constant guidance throughout this project. 
The author would also like to thank Lorenzo Sadun for his kindness in correcting some of the author's mistakes in the background material.
Finally, the author thanks an anonymous referee for pointing out a fatal error in a prior version of this paper in which \autoref{con:nilpotent tiling spaces are fiber bundles} was claimed as a theorem.
I am grateful for their time, energy, and attention to detail.
Their thorough feedback on earlier drafts of this work and their corresponding suggestions have led to substantial improvements in this article.
All remaining confusion and error (fatal or otherwise) belong solely to the author.

\section{Tiling Spaces}\label{sec:tiling spaces}
    In this section, we present the required background on tiling spaces.
The first few subsections are very general, and it is helpful to keep in mind the particular case $G = \mc{M}$ a nilpotent Lie group acting on itself by left translation, with basepoint $m_0 = e$ the origin.
Observe that while we consider left-action, the theory can be made to work for right-actions as well without any substantial change. 
We will eventually focus on well-known results about simple Euclidean tilings.
All tilings are assumed to have finite local complexity (see below) unless otherwise stated.

\subsection{Basic Preliminaries}

    Let $(\mc{M},m_0)$ be a complete, based Riemannian manifold with a distance function $d_{\mc{M}}$.
    A \term{tiling} of $\mc{M}$ is a subdivision of $\mc{M}$ into closed topological disks called \term{tiles} which intersect only along their boundaries.
    We will say that $\mc{T}$ is \term{normal} if tiles are uniformly bounded topological disks which intersect in connected sets.
    We will further say that a tiling is \term{geometrically normal} if, in addition, the intersection of tiles are piecewise smooth submanifolds which are themselves closed topological disks.
    Unless otherwise stated, we assume that all tilings are geometrically normal.
    
    Given a tiling $\mc{T}$ of $\mc{M}$ and a subset $U \sse \mc{M}$, the \term{patch of $\mc{T}$ around $U$} is 
    \[
        [U]_{\mc{T}} = \set{\text{tiles which intersect the closure $\bar{U}$ nontrivially}}.
    \]
    If the tiling $\mc{T}$ is clear from context, we drop the subscript and simply write $[U] = [U]_{\mc{T}}$.
    
    Fix a closed subgroup $G \le \mathrm{Isom}(\mc{M})$ acting transitively on $\mc{M}$. 
    Suppose moreover that $d_\mc{M}$ is left-invariant under $G$, and that $G$ has a norm $\|\cdot\|$ associated to it.
    We will be interested in studying tilings which have nice properties relative to $G$.
    We therefore say that a tiling is a \term{$G$-tiling of $\mc{M}$} when we want to track the action of $G$.
    
    We may often switch between $\mc{T}$ and $\mc{M}$, using each to denote both the manifold as well as the tiling, though we try to use $\mc{T}$ whenever the tiles are relevant.
    We call a set $\mc{P}$ of tiles in $\mc{M}$ a set of \term{prototiles for $\mc{T}$} if, for every tile $T$ in $\mc{T}$, there is a prototile $P \in \mc{P}$ and some $g \in G$ such that $g . P = T$, where $g . P := \set{g.p \sothat p \in P}$.
    Given a tiling $\mc{T}$, we define $g.\mc{T}$ to be the tiling for which all tiles have been translated by $g$. That is, $g.\mc{T} = \set{g.T \sothat T \in \mc{T}}$.

    The action of $G$ defines an equivalence relation on patches of $\mc{T}$. 
    Given patches $[U],[V] \sse \mc{T}$, we say that $[U] \sim_{G} [V]$ if $[U]_{g.\mc{T}} = [V]_{\mc{T}}$ for some $g \in G$.   
    If the group $G$ is clear from context, we will simply write $[U] \sim [V]$.
    
    We say that $\mc{T}$ has \term{finite local complexity} (or simply \term{FLC}) if for any compact $K \sse \mc{M}$, the set of equivalence classes of translates of $K$ is finite. That is,
    \[
        \# \set{[g.K] \sothat g \in G}/\sim \; \; < \infty
    \]
    for every compact set $K$.
    We will often consider compact sets of the form $K = B_r(m_0)$ where $m_0 \in \mc{M}$ is the basepoint.
    As $G$ is a transitive action, and $d_\mc{M}$ is left-invariant, $\mc{T}$ has FLC precisely when
    \[
        \# \set{[B_r(x)] \sothat x \in \mc{M}}/\sim \; \; < \infty
    \]
    for all $r > 0$.
    Observe for example that if $\mc{T}$ has FLC, then $\# \mc{P} < \infty$. 
    All tilings throughout the paper are assumed to have finite local complexity.

    Here, we define a metric, called the \term{tiling metric}, on all $G$-tilings of $\mc{M}$.
    We want this metric to describe when one tiling looks like another at sufficiently large scales around the basepoint, up to a "small" transformation by $G$.
    If $\mc{T}$ and $\mc{T}'$ are $G$ tilings of $\mc{M}$, we say that they are \term{$r$-agreeable} if there exist $g,g' \in G$ with $\|g\|,\|g'\| < \frac{1}{2r}$ such that $[B_r(m_0)]_{g.\mc{T}} = [B_r(m_0)]_{g'.\mc{T}'}$. 
    In other words, the two tilings agree on a patch of radius $r$ around the basepoint, up to a relatively small translation.
    We define
    \[
        R(\mc{T},\mc{T}') := \sup \set{r > 0 \sothat \mc{T} \text{ and } \mc{T}' \text{ are $r$-agreeable}}.
    \]
    The tiling metric is then defined by
    \[
        d_{G,\mc{M}}(\mc{T},\mc{T}') := \min\set{1,\frac{1}{R(\mc{T},\mc{T}')}}.
    \]
    A \term{tiling space} is a collection $\ms{T}$ of tilings which is complete with respect to the tiling metric, and invariant under the action of $G$, so that $g.\mc{T} \in \ms{T}$ for every $\mc{T} \in \ms{T}$.
    Given a tiling $\mc{T}$, the smallest tiling space containing $\mc{T}$ is called the \term{hull} of $\mc{T}$, and will be denoted by $\Hull{G,\mc{T}}$, or simply $\Hull{\mc{T}}$ if context allows.
    Equivalently, if $\mc{O}(\mc{T}) := \set{g.\mc{T} \sothat g \in G}$ is the orbit of $\mc{T}$ under the action of $G$, then $\Hull{G,\mc{T}} = \bar{\mc{O}(\mc{T})}$.
    Though we are interested in tiling spaces in general here, it is often instructive to keep the hull of a particular tiling in mind as a primary example of such spaces.

\subsection{G\"{a}hler Complexes and Inverse Limits}

    The following result from \cite{sadun2003tilinginverselimits} shows that a large class of tiling spaces can be also be constructed as the inverse limit of certain branched manifolds.

    \begin{theorem}[{\cite[\textsection 2]{sadun2003tilinginverselimits}}]\label{thm:tiling spaces are inverse limits}
        Let $\ms{T}$ be a geometrically normal, FLC $G$-tiling space over a based Riemannian homogeneous space $\mc{M}$. 
        Then for each $n \in \NN$ there is a branched manifold $\Gamma_{n}$ with $\dim \Gamma_{n} = \dim G$, and a continuous map $\sigma_{n} : \Gamma_{n} \to \Gamma_{n-1}$, such that $\ms{T} = \invlim (\Gamma_{n},\sigma_n)$.
    \end{theorem}

    \begin{note}
        The theorem as stated in the original paper does not make any reference to the normality of the tilings. 
        We add in the assumption that the tiling space is \term{geometrically normal} for reasons to be explained in \autoref{ssec:simple tilings}.
    \end{note}

    We do not reproduce the proof here, but we mention briefly how these branched manifolds $\Gamma_{n}$, known as \term{G\"{a}hler complexes}, are constructed. 
    Their role is to encode in a single object the combinatoric structure of the set of iterative layers of collars around prototiles $\mc{P}$. 
    Taking $n \to \infty$ represents looking at the combinatorics of larger and larger snapshots around a particular tile in $\mc{T}$. 
    These complexes are defined by setting an equivalence relation $\simeq_n$ on a tiling $\mc{T}$ by $x \simeq_n y$ whenever the "$n$-coronae" around $x$ and $y$ are identical after superimposing $x$ over $y$ by the action of $G$.

    Here we provide a bit more detail, with the specific goal of understanding the construction of $\Gamma_1$, since this object will be used frequently in following sections. 
    Fix a tiling $\mc{T}$ of $\mc{M}$, and let $x \in \mc{M}$. 
    Set $[x]_{0} := \set{T \sse \mc{T} \sothat x \in T}$. 
    This is the patch of all tiles in $\mc{T}$ which contain $x$. 
    For example, if $x$ belongs to a boundary common to two tiles $T_1$ and $T_2$, then $[x]_{0}$ is the union $T_1 \cup T_2$ glued along the common face containing $x$. 
    If $g \in G$ is such that $g.x = y$ and such that the patches of tiles $g.[x]_{0}$ and $[y]_0$ are identical, we will abuse notation to say $[x]_0 = [y]_0$.

    Inductively, having defined $[x]_{n-1}$, let $[x]_n := \set{T \sothat T \cap [x]_{n-1} \neq \emptyset}$. This is called the \term{$n$th corona} around $x$, and contains the patch of tiles that nontrivially intersect the $(n-1)$-corona around $x$. 
    Again, we write $[x]_{n} = [y]_{n}$ if some action of $G$ aligns both the basepoints $x$ and $y$, and the $n$-coronae around them.
    
    More generally, we can define $[U]_n$ for any subset $U$ in a similar way. 
    If $g.[U]_n = [V]_n$ for some subset $V$ and some $g \in G$, then we will write $[U]_n = [V]_n$. 
    When $U = T^{int}$ is the interior of a tile, we call the class $[T^{int}]_n$ an \term{$n$-collared prototile}.
    When $n = 1$, we simply call the class $[U]_1$ a \term{collared prototile}. 
    Given a specific tiling $\mc{T}$ we define $\mc{P}_n := \set{[U]_n \sothat U \sse \mc{T} \text{ is the interior of some tile}}$ to be the set of \term{$n$-collared prototiles}. 
    For example, there is an obvious bijection $\mc{P}_0 \onetoone \mc{P}$.
    In fact, when speaking of an $n$-collared prototile $[T^{int}]_n$, we will often want to refer to the geometry of the underlying (uncollared) base prototile, by which we mean the prototile $P \in \mc{P}$ representing $[T^{int}]_0 \in \mc{P}_0$ under this correspondence.

    For each $n$, then, we have an equivalence $\sim_n$ on tiles of $\mc{T}$ and we define $\Gamma_{n} := \mc{T}/\sim_n$. 
    These are essentially the branched manifolds of the above theorem.
    A point $[x]_n \in \Gamma_{n}$ describes the $n$th corona around $x \in \mc{T}$. 
    When $n = 1$, this gives a description of the base prototile to which a point belongs, as well as the shapes of the (uncollared) tiles surrounding the base prototile. 
    
    \begin{definition}
        The complex $\Gamma_n$ of a tiling space $\ms{T}$ is the $n$th \term{G\"{a}hler complex} of $\ms{T}$.
        When $n = 1$ the complex $\Gamma = \Gamma_1$ is known as the (once-collared) \term{Anderson-Putnam complex}. 
    \end{definition}

    As long as our tiling is geometrically normal, the prototiles may only intersect in piecewise smooth topological disks. 
    As a tiling with FLC, there are only finitely many types of intersections, and so each $\Gamma_{n}$ is a cellular complex where each cell corresponds precisely to a $n$-collared prototile with appropriate cellular subdivision along the boundary.

    Observe that this construction extends to tiling spaces rather than to just a specific tiling. Indeed, the $n$-collared prototiles of a tiling space $\ms{T}$ are the collection of $n$-collared prototiles that show up in some tiling $\mc{T}$ in $\ms{T}$. 
    Two $n$-coronae will be equivalent if they are translates of each other in $\mc{M}$ after aligning basepoints by some element $g \in G$. 
    The notion of $n$-collared prototiles extends in the same manner.

    The branched manifold $\Gamma$ comes equipped with local charts $U \to \Gamma$ where $U = \set{B_i}$ can be taken to be a (finite) collection of metric balls $B_i \subseteq \mc{M}$ called \term{coordinate sheets}.
    For small enough balls the metric on $\mc{M}$ induces a diffeomorphism from each $B_i$ to a Euclidean disk.
    
\subsection{Notions of Equivalence}
    As topological objects, tiling spaces may be considered equivalent if they are homeomorphic. 
    This is naturally a very weak notion of equivalence, though, and loses out on the dynamical properties of the acting group $G$, as well as the specific tiling on which we act.

    Two tiling spaces $\ms{T}_1$ and $\ms{T}_2$ will be called \term{topologically conjugate} if there is a homeomorphism $\phi : \ms{T}_{1} \to \ms{T}_{2}$ such that $g.\phi(\mc{T}) = \phi(g.\mc{T})$ for any $\mc{T} \in \ms{T}_1$.
    Though this provides some more structure on the dynamics of our tiling space, we have still lost out on the actual tilings associated to them.

    A tiling $\mc{T}'$ will be called \term{locally derivable from $\mc{T}$} if, for some $r > 0$, whenever $[B_r(x)]_{g.\mc{T}} = [B_r(y)]_{\mc{T}}$ for any $x,y \in \mc{M}$ and $g \in G$, then $[B_1(x)]_{g.\mc{T}'} = [B_1(y)]_{\mc{T}'}$. 
    That is, whenever two patches in $\mc{T}$ of large enough size agree after a translation by $G$, so will the local patches in $\mc{T}'$.
    If $\mc{T}$ and $\mc{T}'$ are each locally derivable from one another, they are called \term{mutually locally derivable}.

    This notion extends to tiling spaces as well. 
    Suppose that $\phi : \ms{T}_1 \to \ms{T}_2$ is a topological conjugacy.
    We say that $\ms{T}_2$ is \term{locally derivable from $\ms{T}_1$ via $\phi$} if there is a radius $r > 0$ such that for all $\mc{T},\mc{T}' \in \ms{T}_1$, whenever $[B_r(m_0)]_{\mc{T}} = [B_r(m_0)]_{\mc{T}'}$, then $[B_1(m_0)]_{\phi(\mc{T})} = [B_1(m_0)]_{\phi(\mc{T}')}$.
    If $\ms{T}_1$ is locally derviable from $\ms{T}_2$ via $\phi$ as well, we say that $\ms{T}_1$ and $\ms{T}_2$ are \term{mutually locally derivable}, or simply \term{MLD} (via $\phi$).

\subsection{Simple Tilings in \texorpdfstring{$\RR^n$}{Euclidean Space}}\label{ssec:simple tilings}
    In this subsection, we restrict our attention explicitly to the Euclidean setting $G = \mc{M} = \RR^n$ acting on itself by translation.

    A \term{simple tiling} is an FLC tiling of $\RR^n$ by convex polytopes meeting full-face to full-face. 
    Often this is described more loosely as tilings by polyhedra meeting full-face to full-face. 
    We add the convexity assumption as a useful (but technically unnecessary) condition.
    The assumption that a tiling is simple is widely used to help control both the geometry and combinatorics of otherwise unwieldy tile types.
    It is well-known that such a restriction does not lose us any information in $\RR^n$.

    \begin{theorem}[{\cite[\textsection 4.2]{priebe2000towards}}]\label{thm:normal tilings are simple}
        Every normal tiling of $\RR^n$ is MLD to a tiling by Vorono\"{i} cells, and therefore is MLD to a simple tiling. 
        \begin{sketch}
            We omit the full proof of this fact here, noting that there will be similarities when we prove \autoref{thm:non-Euclidean normal tilings are simple} later one.
            The key ingredients are that every tile can be represented by a "locator point" (such as the barycenter of the tile). 
            The collection of these points form a Delaunay set $D_{\mc{T}}$.
            The Derived Vorono\"{i} tiling associated to $D_{\mc{T}}$ is MLD to $\mc{T}$ by construction, and since Vorono\"{i} cells in $\RR^n$ are convex polytopes meeting full-face to full-face, such a tiling is simple.
        \end{sketch}
    \end{theorem}

    Observe that the definition of a simple tiling does not immediately generalize to the setting of a general $G$-tiling of a manifold $\mc{M}$.
    In particular, convex polytopes may not be well-behaved at large scales.
    In fact, it is not even clear that Vorono\"{i} cells may be deemed "polytopes" in any reasonable sense.
    For more on this, we invite the reader to consider \cite{boissonnat2016obstruction} wherein it is shown that the Vorono\"{i} faces of a Vorono\"{i} diagram may not even form closed topological disks.
    The example manifold provided is non-homogeneous, but can be made to approximate Euclidean space arbitrarily well.

    \begin{note}
        It is for this reason that we hesitate to take \autoref{thm:tiling spaces are inverse limits} in the generality as originally stated, prompting us to add the extra "geometric normality" assumption. 
        This assumption is essentially assumed in first half of the proof, but is dropped later on, claiming that the Vorono\"{i} cells of a locator set will induce a "polytope" cell division.
        From this it is implied that the tiles are closed topological disks meeting in connected, topological disks after sufficient subdivision.
        
        As explained above, this is not guaranteed for Vorono\"{i} cells in a general Riemannian manifold.
        We admit that perhaps the assumption that $\mc{M}$ is a Riemannian \term{homogeneous} space will allow for a "polytope" Vorono\"{i} cell decomposition, but the author has been unable to locate a proof of this fact in the non-Euclidean setting. 
        If the claim is indeed true in the general setting, then we have not added any unnecessary assumptions.
    \end{note}

    Using simple tilings allows us to focus on how tiles are allowed to meet, rather than on the specific geometry of their intersecting boundaries.
    Because of this, simple tiling spaces are often easier to understand, and the assumption that tiling spaces are simple is practically a given in the field of study.
    For example, the topology of such tiling spaces is fairly well-understood (see \cite{sadun2008topology}). 
    As a specific case-study, a theorem of \cite{sadun2003tilingfiberbundles} states that such tiling spaces are fiber bundles over the torus $\TT^n$.
    We demonstrate progress towards a generalization of this theorem, using our definition of a simple tiling, in \autoref{sec:fiber bundles}.

\subsection{Examples}

Here we mention a few classic examples of tilings of $\RR^2$, along with some of the remarkable features of each. 
With the exception of the "hat" tiling, exposition on these examples can be found in \cite{sadun2008topology}, and we refer the interested reader to this book for descriptions of these tilings.

    \begin{example} 
    \begin{enumerate}[before=\leavevmode, beginpenalty=10000]
        \item In the \term{chair tiling}, there are $4$ distinct prototiles: one for each orientation of the "chair". These tiles are polygonal, but do not technically meet full-face to full-face, and are not even convex. This can be modified either by replacing chairs by "arrow" tiles, or by an appropriate subdivision of (the faces of) the tiles. 
        \item A tiling of $\RR^2$ by squares of unit side length above the $x$-axis, and by squares of irrational dimensions (say, side length $\sqrt{2}$) below the $x$-axis, will exhibit a "fault line" along this axis. Such a tiling cannot be rectified by appropriate subdivision of the faces of tiles without introducing infinitely many prototiles. Since there are infinitely many ways for tiles to meet, this tiling has \term{infinite local complexity (ILC)}.
        \item A classical \term{Penrose tiling} of $\RR^2$ can be created by two rhombi, each in five different orientations, yielding a total of $10$ tile types. Due to local matching rules however, one often distinguishes a rhombus from its rotation by $180$ degrees. Under this convention there are $20$ prototiles.
        \item A variation on a Penrose tiling of $\RR^2$ uses chicken-shaped tiles in place of rhombi. Though the chicken-tiling is non-polygonal, it is still a normal tiling, and therefore MLD to some simple tiling. In fact, it is MLD to a Penrose tiling. 
        \item Though the \term{pinwheel tiling} has only one prototile up to isometry, this one prototile shows up in infinitely many orientations. Therefore, up to translation, this tiling has infinitely many prototiles.
        Extending the group $G$ acting on $\RR^2$ to be the group of Euclidean motions, there would be only one prototile.
        \item The recently discovered triskaidecagonak einstein (the $13$-sided polygonal tile from \cite{smith2023aperiodic}) known as "the hat" can only tile $\RR^2$ aperiodically. This tile has finitely many admissible orientations (including rotations and reflections) in a such a tiling, and therefore any such tiling has finitely many prototiles.
        In order for the hat to meet the "full-face to full-face" condition of a simple tiling, one must subdivide one of the edges to create a $14$-gon, similar to the subdivision of the chair tiling.
        Although the convexity assumption is not satisfied, we still consider the resulting tiling simple.
        \item A \term{periodic tiling} $\mc{T}$ of $\RR^n$ is one for which there is a cocompact lattice $\Lambda \le \RR^{2}$ such that $\lambda \cdot \mc{T} = \mc{T}$ for all $\lambda \in \Lambda$. 
        For example, the tiling of $\RR^{2}$ by squares whose vertices lie on integer points is periodic. 
        This definition extends in the natural way to any $G$-tiling of a manifold $M$.
        \qedhere{}
    \end{enumerate}
    \end{example}

    We conclude this section by offering a few examples of mathematical quasicrystals in the connected, simply connected, rational, nilpotent Heisenberg group $\HH$. 
    In particular, $\HH$ should not be confused with hyperbolic space or the quaternions.
    Some details of such groups may be found in \autoref{sec:nilpotent lie groups}. 
    
    We note that the quasicrystals described below are most effectively described as Delaunay sets rather than as tilings.
    The interested readers may take it upon themselves to construct specific prototiles associated to such Delaunay sets.
    
    \begin{example}
        The Heisenberg group $\HH$ can be naturally identified with upper triangular $3 \times 3$ matrices with $1$s along the diagonal, with entries in $\RR$. 
        There is a cocompact lattice $\HH(\ZZ)$ identified with those matrices of integer values.
        
        A periodic simplicial tiling in $\HH$ by a fundamental domain exists, and the vertex set of the tiling is $\HH(\ZZ)$.

        The hull of this tiling will be homeomorphic to the Heisenberg manifold $\Heis = \HH/\HH(\ZZ)$. 
        More generally, if $\mc{T}$ is a periodic tiling of a simply connected, connected, rational, nilpotent lie group $\GG$ with $e$ a vertex, then the vertex set of $\mc{T}$ will be a lattice $\Lambda$ inside $\GG$, and it can be seen that $\Hull{\mc{T}} = \GG/\Lambda$ is a nilmanifold.
    \end{example}

    A more complicated example of an FLC Delaunay set can be found in \cite[Ex. 4.18]{kaiser2022complexity}, which we reproduce here.
    Again, though there is not necessarily one particular tiling associated to this Delaunay set, by taking (some approximation to) Vorono\"{i} cells around each point, we might hope to produce prototiles for a tiling, or an actual tiling itself.
    We warn that such a "psuedo-Vorono\"{i}" tiling will depend on the particular Riemannian metric given to $\HH$.
    Nonetheless, this example provides a nontrivial FLC Delaunay set which may be instructive to keep in mind.
    For full details, we direct the reader to the original paper, wherein the cut-and-project method is fully described.

    \begin{example}\label{ex:cut-and-project Heisenberg}
        Consider the product of Heisenberg groups $\HH \times \HH$.
        Let \[L = \set{(a,b,c,a^*,b^*,c^*) \sothat a,b,c \in \ZZ[\sqrt{2}]} \sse \HH \times \HH,\] where $x^*$ denotes the Galois conjugate of $x$. 
        Let $W = \exp(I_{1/2} \times I_{1/2} \times I_{1/2})$, where we view $I_r = [-r,r]$ as the image of the interval in $\mf{h} = \RR^3$ under the exponential map.
        Using a cut-and-project scheme with this window $W$, we arrive at an FLC Delaunay set $D$ within $\HH$.

        In $D$, the points are aligned along the $x$, $y$, and $z$ directions in "straight" lines (when we consider them as lines in $\RR^3$).
        For example, given a point $(a,b,c) \in D$ with $a = \alpha_1 + \alpha_2\sqrt{2}$. If a point $(a',b,c) \in D$ is close to $(a,b,c)$, then (because the only change is in the first coordinate) we will either have $a' = a \pm (1 + \sqrt{2})$ or $a' = a \pm (2 + \sqrt{2})$.
        It can be seen that up to translation there are only finitely many patches of any given radius, and that balls of sufficiently large radius cover all of $\HH$ when translated to points of $D$.
    \end{example}

\section{Delaunay Sets and Triangulations by Riemannian Simplices}\label{sec:delaunay}
    We now develop some technical background for the proof of \autoref{thm:non-Euclidean normal tilings are simple}. 
The idea of the proof is to develop a pattern-equivariant ("PE") indicator set for a tiling space $\ms{T}$, which in turn can be used to create a PE triangulation of any tiling $\mc{T} \in \ms{T}$.
This process essentially boils down to identifying a Delaunay set $\mc{D}$ on the (once collared) Anderson-Putnam complex $\Gamma$ which will in turn induce a triangulation of $\Gamma$ (or, rather, of coordinate patches of $\Gamma$).
Moreover, we would like $\mc{D}$ itself to serve as the set of vertices of the triangulation, so that one can discuss either the vertices of a simplex (the point pattern) or the simplex itself (the geometric and combinatorial behavior) without any loss of information.
This will imply that the resulting triangulation will be MLD to the one which determined $\mc{D}$ in the first place.

Our proof of \autoref{thm:non-Euclidean normal tilings are simple} technically only requires that the indicator set be generic in the sense of Delaunay; that is, the points are in general position so as to induce a natural Delaunay triangulation.
However, by adding greater control over certain parameters of this genericity, we will be able to further develop PE triangulations which remain stable under sufficiently small PE perturbations.
This feature of perturbation, while irrelevant for \autoref{thm:PE cohomology}, is a crucial step in our progress towards \autoref{con:nilpotent tiling spaces are fiber bundles}.

Building on the work of \cite{dyer2015riemannian}, an algorithm is developed by Boissannat, Dyer, and Ghosh in \cite{boissonnat2013stability,boissonnat2018delaunay} for developing such a Delaunay set in the context of compact Riemannian manifolds (more specifically, see \cite[Theorem 3]{boissonnat2018delaunay}).
One benefit of this approach is that the indicator set and the triangulation uniquely define one another.
For example, each simplex is a \term{Riemannian simplex}, which is determined precisely by the location of its vertices, and does not rely on information about neighboring points.

Since the input to the algorithm is simply a finite atlas for the manifold, it is not too difficult to see that the algorithm may be modified to produce a Delaunay set on $\Gamma$ whose pullback to $\mc{M}$ induces a triangulation in the same manner.
The main difficulty comes by ensuring that the selected Delaunay set on $\Gamma$ induces a triangulation on every sheet around a branch set simultaneously; that is, the Delaunay set must induce a triangulation which is in some sense consistent across branches of $\Gamma$. 
This is essentially accomplished by considering thickened branch sets (which by assumption are initially piecewise smooth).\footnote{In essence, this is akin to considering a Barge-Diamond complex corresponding to $\mc{T}$, rather than the Anderson-Putnam complex.}

With this in mind, the rest of this section is dedicated to setting up the precise result of the algorithm of Boissannat, Dyer, and Ghosh, together with some of the necessary vocabulary and stability results which will allow us to prove \autoref{thm:non-Euclidean normal tilings are simple}.

\subsection{Riemannian Simplices}

In $\RR^{n}$, a generic set of $n+1$ points uniquely determines a simplex.
Naturally one might consider a simplex to be the convex hull of its vertices.
The standard Euclidean $n$-simplex $\eusimplex^{n}$ may be defined via \term{barycentric coordinates}:
$${\displaystyle \eusimplex^{n}=\left\{(t_{0},\dots ,t_{n})\in \mathbb {R} ^{n+1}~{\Bigg |}~\sum _{i=0}^{n}t_{i}=1{\text{ and }}t_{i}\geq 0{\text{ for }}i=0,\ldots ,n\right\}}.$$

The definition via barycenters may be generalized beyond the Euclidean setting.
Let $\mc{M} = \mc{M}^{n}$ be a manifold of bounded geometry with distance function $d_\mc{M}$. 
In particular, let $\iota$ denote the injectivity radius of $\mc{M}$ and let $\kappa$ denote the sectional curvature, and suppose $\iota \ge c > 0$ is bounded away from $0$ and that $|\kappa| \le C < \infty$ is absolutely bounded above.

\term{Riemannian simplices} are introduced in \cite{dyer2015riemannian} as simplices which will be determined by their vertex set in $\mc{M}$ under suitable conditions.
Take a finite set $\sigma^{j} = \set{p_0,p_1,\dots,p_j}$ contained in an open geodesic ball $B_{r} \subseteq \mc{M}$ of radius $r > 0$ for which $\bar{B_{r}}$ is convex.
That is, for any $x,y \in B_{r}$ there is a minimizing geodesic $\gamma(t)$ from $x$ to $y$ which is unique in $\mc{M}$, and moreover $\gamma(t) \in B$ for all $t$.
Let $\lambda = (\lambda_0,\lambda_1,\dots,\lambda_n) \in \eusimplex^{n}$. 
Define $\mc{E}_{\lambda} : \bar{B}_{r} \to \RR$ by 
\begin{equation}
    \mc{E}_{\lambda}(x) = \frac{1}{2}\sum_{i} \lambda_i d_M(x,p_i)^{2}.
\end{equation}
The \term{barycentric center of mass} is $\min_{x \in B_{r}}\mc{E}_{\lambda}(x)$. 
There is a constant $C(\iota,\kappa)$ such that for $r < C(\iota,\kappa)$, the barycentric center of mass exists and is unique, as shown in \cite[Thm. 1.2]{karcher1977riemannian} (see also \cite[Lemma 3]{dyer2015riemannian}).
The image of the map 
\begin{align*} \mc{B}_{\sigma^{j}} 
    : \eusimplex^j & \to \mc{M} \\ 
    \lambda & \mapsto \operatorname{argmin}_{x \in \bar{B}_{r}} \mc{E}_{\lambda}(x) 
\end{align*} 
is denoted $\rsimplex^{j}_{\mc{M}}$ and is called a \term{Riemannian simplex} whose vertices are $\sigma^j$. 
The map $\mc{B}_{\sigma^j}$ is called the \term{barycentric coordinate map} of $\rsimplex^{j}_{\mc{M}}$.
If $\mc{B}_{\sigma^j}$ is a smooth embedding, we say that $\rsimplex^{j}_{\mc{M}}$ is \term{nondegenerate}; otherwise we say that $\rsimplex^{j}_{\mc{M}}$ is \term{degenerate}.
We may also speak of the configuration of vertices $\sigma^{j}$ itself as being (non)degenerate if the corresponding Riemannian simplex is.
An \term{$i$-face} of $\rsimplex^{j}_{\mc{M}}$ is the image of an $i$-face of $\eusimplex^{j}$.

The upshot of all of this is that under sufficiently nice geometric restrictions on $\mc{M}$, a set of $n+1$ points in a small enough ball always uniquely induce a (possibly degenerate) Riemannian simplex.
The work of \cite{dyer2015riemannian} provides further genericity conditions on $\sigma^{j}$ under which the corresponding Riemannian simplex is guaranteed to be nondegenerate. 
In this case we will often identify the simplex $\rsimplex^{j}_{\mc{M}}$ with the vertex set $\sigma^{j}$.

\subsection{The Delaunay Complex and Stability of Triangulations}
A key contribution of \cite{boissonnat2013stability} to this story is the application of the above ideas to triangulating $\mc{M}$ from a given Delaunay set $\mc{D} \subseteq \mc{M}$.
Among the subtlety addressed is the fact that not every Delaunay set---even with a sufficiently dense sample---will induce such a triangulation (in contrast to claims such as those of \cite{leibon2000delaunay}).
In order to distinguish between subsets of $\RR^{n}$, and of $\mc{M}$ more generally, we will use bold face font to denotes sets of $\RR^{n}$, and non-bold font for a more general manifold.
Sufficient conditions on a Euclidean Delaunay set $\mathbf{D} \subseteq \RR^{n}$ are given in \cite{boissonnat2013stability} for which a triangulation is guaranteed to exist.
This thread is further picked up in \cite{boissonnat2018delaunay} and extended to Delaunay sets in compact Riemannian manifolds, yielding an algorithm for perturbing a given Delaunay set into sufficiently generic position from which such a triangulation may be induced. 
In fact, the simplices of the resulting triangulation in $\mc{M}$ are nondegenerate Riemannian simplices as a consequence of the genericity conditions which we now develop. 

\begin{definition}
    Let $Y$ be a metric space.
    A closed subset $X \sse Y$, is said to be \term{$r$-discrete} if $d(x,x') \ge r$ for all $x \neq x'$ in $X$. 
    If $X$ is $r$-discrete for some unspecified $r$, we will simply say that $X$ is \term{uniformly discrete}.
    On the other hand, $X$ is $R$-dense if $\dist(x,X) \le R$ for all $x \in Y$. 
    If $X$ is $R$-dense for some unspecified $R$, we will simply say that $X$ is \term{coarsely dense}.
    If $X$ is uniformly discrete and coarsely dense, we call $X$ a \term{Delaunay set}. 
    We will also say that $X$ is a \term{$(\mu,\xi)$-net} or that $X$ has \term{mesh $(\mu,\xi)$} if $X$ is a Delaunay set which is $(\mu\xi)$-discrete and $\xi$-dense.
\end{definition}
In the event that $\mathbf{X} \subseteq \RR^{n}$, we introduce finite analogues of these definitions. 
\begin{definition}
    We call $\mathbf{X}$ an $\varepsilon$-\term{sample for a bounded subset} $B \subseteq \RR^{n}$ if $d_{\RR^{n}}(x,\mathbf{X}) < \varepsilon$ for all $x \in \bar{B}$. 
    We also say simply that $\mathbf{X}$ is an $\varepsilon$-\term{sample} if $d_{\RR^{n}}(x,\mathbf{X} \cup \boundary \operatorname{Conv}(\mathbf{X})) < \varepsilon$ for all $x \in \operatorname{Conv}(\mathbf{X})$.
    We may say that $\mathbf{X}$ is a \term{sample set} without reference to $\varepsilon$ if $\mathbf{X}$ is a $\varepsilon$-sample set for some $\varepsilon > 0$.
    (Note that $\operatorname{Conv}(\mathbf{X})$ denotes the convex hull of $\mathbf{X}$.)
    Further, $\mathbf{X}$ is $\lambda$-\term{separated} if $d_{\RR^{n}}(x,y) > \lambda$ for all $x,y \in \mathbf{X}$. 
\end{definition}

In order to turn a sample set into a simplicial complex, we follow the ideas of \cite{boissonnat2013stability}, from which we borrow our notation and terminology.
We will primarily require language about finite "patches" of Delaunay sets in local coordinate systems of $\mc{M}$.
In what follows, we let $\mathbf{P} \subseteq \RR^{n}$ be an $\varepsilon$-sample set; it is convenient to think of this as a coordinate patch of a Delaunay set in $\mc{M}$.

\begin{definition}
    An \term{empty ball} (relative to $\mathbf{P}$) is a ball $B \subseteq \RR^{n}$ which contains no points of $\mathbf{P}$, and a \term{Delaunay ball} is a maximal empty ball. That is, $B = B_{\RR^{n}}(x;r)$ is a Delaunay ball if any empty ball centered at $x$ is contained in $B$. 
    A simplex $\sigma \subseteq \RR^{n}$ is a \term{Delaunay simplex} if there is some Delaunay ball $B$ such that the vertices of $\sigma$ belong to $\mathbf{P}$ and lie on the boundary of the closed ball. 
    That is, $\sigma^{(0)} \subseteq P \cap \boundary \bar{B}$.
    
    The (abstract) \term{Delaunay complex} for $\mathbf{P}$ is the abstract simplicial complex $\Del(\mathbf{P})$ whose simplices are the Delaunay simplices of $\mathbf{P}$.
    We denote by $|\Del(\mathbf{P})|$ the geometric realization of $\Del(\mathbf{P})$.
    We say that $\Del(\mathbf{P})$ \term{embeds} into $\RR^{n}$ if the (continuous) barycentric coordinate map $|\Del(\mathbf{P})| \to \RR^{n}$ which extends the inclusion $\mathbf{P} \into \RR^{n}$ is injective.
    In this case, we identify $|\Del(\mathbf{P})|$ with its image in $\RR^{n}$, and we call $|\Del(\mathbf{P})|$ the \term{Delaunay triangulation} relative to $\mathbf{P}$.

    These definitions extend in the obvious ways to any metric $d$ on $\RR^{n}$, and under such a metric we define $\Del_d(\mathbf{P})$ to be the corresponding Delaunay complex, and (should it exist) $|\Del_{d}(\mathbf{P})|$ the Delaunay triangulation relative to $\mathbf{P}$.
\end{definition}

It is a seminal result of Delaunay that Delaunay triangulations relative to $\mathbf{P} \subseteq \RR^{n}$ exist as long as there is no empty ball $B$ with $n + 2 \le  \# (\mathbf{P} \cap \boundary \bar{B})$; that is, the closure of an empty ball contains at most $n+1$ points of $\mathbf{P}$.
Such restrictions are not too strict in $\RR^{n}$; an arbitrarily small affine perturbation can transform a given Delaunay set into one that is generic, as proved in \cite{delaunay1934sphere}.
This notion of genericity is parametrized in \cite{boissonnat2013stability} as follows:

\begin{definition}
    A Delaunay simplex $\sigma$ with Delaunay ball $B$ is $\delta$-\term{protected} for some $\delta \ge 0$ if $d_{\RR^m}(q, \boundary B) > \delta$ for all $q \in \mathbf{P} \setminus \sigma$.
    In this case, we call $B$ a \term{$\delta$-protected Delaunay ball} for $\sigma$. 
    If the affine hull of $\mathbf{P}$ is all of $\RR^{n}$ and every Delaunay $n$-simplex is $\delta$-protected for some $\delta \ge 0$, we say that $\mathbf{P}$ is $\delta$-\term{generic}. 
    In the event the parameter $\delta \ge 0$ is unspecified, we will simply call $\mathbf{P}$ \term{generic}.
\end{definition}

Note that we allow $\delta = 0$ in the above definition. This corresponds to Delaunay's original definition of genericity.
However in this case, a $0$-protected simplex $\sigma$ may be "fragile" in the sense that even arbitrarily small perturbations of $\mathbf{P}$ may cause the combinatorics of the Delaunay triangulation around $\sigma^{(0)}$ to change drastically.

By contrast, the payoff of finding a simplex $\sigma$ with a $\delta$-protected ball for $\delta > 0$ becomes self-evident: sufficiently small perturbations preserve the combinatorics of the Delaunay triangulation of $\mathbf{P}$.
In other words, in the case $\delta > 0$ a Delaunay triangulation ought to remain "stable" under slight changes either to the metric or to the points themselves.
This is one of the main contributions of \cite{boissonnat2013stability}.
The idea of such stability will be relevant in the proof both of \autoref{thm:non-Euclidean normal tilings are simple}, and is key in making progress towards \autoref{con:nilpotent tiling spaces are fiber bundles}.
We include some of the deatails here for completeness.

\begin{definition}
    Let $\mathbf{P} \subseteq \RR^{n}$ be a (possibly finite) $\varepsilon$-sample set. The set of \term{deep interior points} is 
    \[
        \mathbf{P}_{I} = \set{p \in \mathbf{P} \sothat d_{\RR^{n}}(p,\boundary \operatorname{Conv}(P)) \ge 4\varepsilon}.
    \]
    We say that $\mathbf{P}$ is $\delta$-\term{generic for a subset} $\mathbf{P}_{J} \subseteq \mathbf{P}_{I}$ if every $n$-simplex of $\star(\star(\mathbf{P}_{J} ;\Del(\mathbf{P})))$ is $\delta$-protected. 
    The \term{safe interior simplices} are the simplices of $\star(\mathbf{P}_{J} ;\Del(\mathbf{P}))$.
\end{definition}

For example, if $\mathbf{P}$ is a $\varepsilon$-sample of a Delaunay set which is globally $\delta$-generic, and $\operatorname{Conv}(\mathbf{P})$ is sufficiently large relative to $\varepsilon$ (e.g., $B_{5\varepsilon}(0) \subseteq \operatorname{Conv}(\mathbf{P})$), then $\mathbf{P}$ is $\delta$-generic for a singleton near $0$.
The ultimate utility of this finite version of $\delta$-genericity is the following stability condition:

\begin{theorem}[{\cite[Theorem 4.14 + Lemma 4.10]{boissonnat2013stability}}]\label{thm:stable stars of point perturbation}
    Let $\mathbf{Q} \subseteq \mathbf{P}$ be a set of interior points, where $\mathbf{P}$ is a $\varepsilon$-sample set, and where $\mathbf{P}$ is $\delta$-generic for $\mathbf{Q}$ for a sufficiently small $\delta$ relative to the mesh of $\mathbf{P}$. Suppose moreover every $m$-simplex of $\mathbf{Q}$ is safe. There is a $\varrho > 0$ (depending on $\varepsilon$) such that if $\zeta: \mathbf{P} \to \mathbf{\tilde{P}}$ is a $\varrho$-perturbation, then 
    \[
        \star(\mathbf{Q};\Del(\mathbf{P})) \stackrel{\zeta}{\cong}  \star(\zeta(\mathbf{Q});\Del(\mathbf{\tilde{P}})).
    \]
\end{theorem}
In other words, we can parametrize the stability of a Delaunay triangulation relative to sufficiently small perturbations. 
The application to general manifolds becomes aparent when we consider a variation on this fact, namely that a Delaunay triangulation is stable under a perturbation on the \textit{metric}. 
\begin{theorem}[{\cite[Corollary 4.18]{boissonnat2013stability}}]\label{thm:stable stars of metric perturbation}
    Let $\mathbf{P}$ be $\delta$-generic for $\mathbf{P}_{J}$ with sampling radius $\varepsilon$, for $\delta$ sufficiently large relative to $\varepsilon$. 
    Moreover, suppose that $d$ is a metric on $U \supseteq \Conv(\mathbf{P})$ and that $|d(x,y) - d_{\RR^{n}}(x,y)| \le \varrho$.
    If $\varrho$ is sufficiently small, then 
    \[
        \star(\mathbf{P}_{J};\Del_d(\mathbf{P})) = \star(\mathbf{P}_{J};\Del(\mathbf{P})).
        \]
    \end{theorem}
In particular, if $\mc{M}$ has bounded geometry, coordinate patches of $\mc{M}$ will have a metric which is quantifiably close to that of $\RR^{n}$, and a sufficiently dense Delaunay set of $\mc{M}$ which is $\delta$-generic (for $\delta > 0$ sufficiently large relative to $\varepsilon$) will have combinatorics identical to those induced by the coordinate patches under the Euclidean metric.
As a consequence, one has stability under perturbations of the points in $\mc{M}$ as well.

\subsection{The Extended Algorithm}

Let $\mc{D} \subseteq \mc{M}$ be a Delaunay set with sample size $\varepsilon > 0$.
We have seen that for $\delta > 0$, there is a Delaunay triangulation $\Del(\mc{D})$ by Riemannian simplices.
Moreover, if $\delta$ is sufficiently large relative to $\varepsilon$, the triangulation is stable under sufficiently small perturbation.

This begs the question: how common is a $\delta$-generic set in $\mc{M}$? 
Such sets can be attained for $\delta = 0$ by a sufficiently small perturbation of $\mc{D} \subseteq \RR^{n}$ as Delaunay showed in \cite{delaunay1934sphere}.
It is the main consequence of \cite{boissonnat2018delaunay} that this program can be carried out more generally in compact manifolds $\mc{M}$ for $\delta > 0$. 
We have rephrased and condensed some of their results slightly for our purposes here.
(We note that while the original definition of a $\delta$-generic set includes $\delta = 0$, the algorithm is designed with $\delta > 0$ in mind---see also \cite{boissonnat2014delaunay}.)

\begin{theorem}[{\cite[Theorem 3]{boissonnat2018delaunay}}]\label{thm:extended algorithm}
    Let $\mc{M} = \mc{M}^{n}$ be a compact manifold with injectivity radius $\Lambda$ and sectional curvature $\kappa$.
    There is a constant $C = C(\Lambda,\kappa,n)$ such that the following holds: 
    If $\mc{D} \subseteq \mc{M}$ is an $\varepsilon$-separated Delaunay set for $\varepsilon < C$, there is a perturbation $\mc{D}'$ of $\mc{D}$ which is $\delta$-generic for some $\delta > 0$.
    Moreover, there is a triangulation $|\Del(\mc{D}')| \to \mc{M}$ by Riemannian simplices. 
    Furthermore, $\Del(\mc{D}')$ is \term{self-Delaunay} (i.e., a Delaunay triangulation of its vertices with respect to the piecewise-linear metric $d_{PL}$ defined by geodesic distances between vertices in $\mc{M}$).
\end{theorem}

In some sense, the result of this algorithm is not all that crucial for our purposes; however, it demonstrates that there is an algorithmic process for perturbing a Delaunay set into one which is $\delta$-generic (for some $\delta > 0$) using coordinate charts from the atlas of $\mc{M}$.
It is this process that we can make use of.
The algorithm is described in detail in \cite{boissonnat2018delaunay}, and we sketch the outline of the process here.

The algorithm begins with a finite atlas, together with transition functions, for a compact manifold $\mc{M}$. We further input a Delaunay set $\mc{D} \subseteq \mc{M}$ whose mesh is sufficiently fine relative to data from the atlas and transition functions, as well as the metric on $\mc{M}$, but otherwise arbitrary.
The set $\mc{D}$ is equipped with an (arbitrary) enumeration.

The algorithm then creates a $\delta$-generic Delaunay set $\mc{D}'$ by visiting each point $p \in \mc{D}$ in turn. The algorithm perturbs $p \mapsto p'$ into a position such that the collection of points which have been perturbed are $\delta$-protected. 
(The existence of a perturbation strategy is demonstrated by a volumetric argument; it is this existence which we most heavily require.)
The algorithm halts when all points of $\mc{D}$ have been perturbed, and the perturbed points $\mc{D}'$ form a $\delta$-generic Delaunay set for $\mc{M}$.

\begin{definition}
    We say that a Delaunay set $\mc{D}$ is \term{good} if the Delaunay simplices of $\mc{D}$ are $\delta$-protected. 
    We say that a point $p \in \mc{D}$ is in \term{good position} relative to some subset $\mc{D}' \subseteq \mc{D}$ if $p \in \mc{D}'$ and $\mc{D}'$ is good.
    We will often not make reference to the specific subset $\mc{D}'$, but the existence of this subset and its relation to $\mc{D}$ should be understood from context.
\end{definition}

Under this notion, we understand that the algorithm proceeds by perturbing each point in turn so that at each stage the collection of perturbed points are good.
Each point is individually perturbed into good position relative to nearby points in a small ball, until all points have been visited.
The algorithm ends when all points have been perturbed into good position.

We note that the order in which points are selected for perturbation is taken as part of the input data, and therefore is something that can be exploited.
In particular, we can choose the order in which to perturb the points of our Delaunay set. 
This feature will allow us to choose a convenient order in which to apply the algorithm to the once-collared Anderson-Putnam complex.

\section{Simple Tilings}\label{sec:simple tilings}
    We now demonstrate how to employ the ideas of \autoref{sec:delaunay} to create an indicator set for a geometrically normal tiling $\mc{T}$ of a (possibly non-compact) manifold $\mc{M}$, which will in turn induce a triangulation of $\mc{M}$ by Riemannian simplices in a PE way.
Although we avoid the technical aspects of the estimates presented in \cite{boissonnat2018delaunay}, we note that the quantities discussed should in fact be algorithmically constructable.
As before, we use $\Gamma$ to denote the once-collared Anderson-Putnam complex of $\mc{T}$.

\subsection{Simple Tilings}

We start by defining simple tilings in a more general context than previously allowed.
To the author's knowledge, such a generalization has not been previously provided.

\begin{definition}
    An FLC $G$-tiling $\mc{T}$ of $\mc{M}$ is said to be \term{simple} if $\mc{T}$ is a self-Delaunay triangulation of $\mc{M}$.
    That is, there is a Delaunay set $\mc{D} \subseteq \mc{M}$ such that $\Del{\mc{D}}$ is well-defined, and such that $\mc{D}$ is the vertex set of the triangulation $|\Del{\mc{D}}| \to \mc{M}$. 
    A tiling space will be called \term{simple} if this condition is true of every $n$-collared prototile for every $n$.
\end{definition}

Clearly any FLC tiling of $\RR^n$ by convex polytopes meeting full-face to full-face will be MLD to a simple tiling in this new definition by appropriately subdividing the polytope into linear simplices.
Therefore any simple tiling in the standard sense for $\RR^{n}$ will be simple in this sense as well.
The main purpose of this section is to prove that an analogue of such a subdivision is possible in the non-Euclidean context.
The following proposition is key in the proof:

\begin{restatable*}{proposition}{genericOnAP}\label{prop:existence of generic across branches}
    One can construct a finite indicator set on $\Gamma$ whose pullback to $\mc{M}$ is a $\delta$-generic Delaunay set of arbitrarily fine mesh for some $\delta > 0$.
\end{restatable*}

Given a finite set $D \subseteq \Gamma$ we will use the phrase \term{coordinate patch} to denote the restriction of $D$ to a coordinate sheet of $\Gamma$ under the metric of $\mc{M}$.
If $B_i \subseteq \mc{M}$ and $\pi_i : B_i \to \Gamma$ is a local coordinate sheet for $\Gamma$, this allows us to view the coordinate patch $D \cap \pi_{i}(B_i)$ both as a subset of $\Gamma$ and of $\mc{M}$ simultaneously.
Because genericity is a local condition, we can abuse terminology slightly to say (a subset of) $D \subset \Gamma$ is $\delta$-generic to mean that all coordinate patches of $D$ (i.e., restrictions of $D$ to sheets of branches) viewed as subsets of $\mc{M}$ are $\delta$-generic.
We may also abuse terminology and call $D$ a Delaunay set on ${\Gamma}$, by which we mean that the coordinate patches of $D$ are Delaunay sets when viewed as subsets of $\mc{M}$.
We use $D$ to denote a finite set on $\Gamma$ and $\mc{D}$ to denote its pullback to $\mc{M}$ whenever relevant; however this potential abuse of notation should be kept in mind.

\simpleTilings

\begin{proof}[Proof (assuming \autoref{prop:existence of generic across branches})]
    The strategy of the proof is not too difficult.
    Take a tiling $\mc{T}$ of $\mc{M}$, and let $\pi : \mc{M} \to \Gamma$ be the natural identification of (collared) tiles to the Anderson-Putnam complex $\Gamma = \Gamma_1$.
    Using \autoref{prop:existence of generic across branches}, take $D$ a finite indicator set on $\Gamma$ such the pullback $\mc{D} := \pi^{-1}(D) \subseteq \mc{M}$ is a $\delta$-generic Delaunay set for some $\delta > 0$.
    Then the Delaunay complex $|\Del{\mc{D}}|$ induced by $\mc{D}$ is well-defined. 

    The complex $\mc{T}_{\mc{D'}} \subseteq \mc{M}$ defined by $|\Del{\mc{D}}| \to \mc{M}$ is a triangulation by Riemannian simplices and evidently will be MLD to $\mc{T}$ because the vertex set $\mc{D}$ simply acts as an indicator set of any given tile by the map $\mc{D} \to D \subseteq \Gamma$.
    In the case of a tiling space more generally, one pulls back the set $D$ to each collared tile individually; these piece together into a $\delta$-generic Delaunay set in any collared tile in any given tiling of the tiling space.
\end{proof}

Note that the above proof demonstrates a slightly stronger claim. We record these details as a separate result, as they will be crucial in \autoref{sec:fiber bundles}.

\begin{proposition}\label{prop:normal tilings are generic and simplicial}
    The simple tiling $\mc{T}_{\Delta}$ produced in \autoref{thm:non-Euclidean normal tilings are simple} is a self-Delaunay triangulation of $\mc{M}$ by a $\delta$-generic Delauay set $\mc{D} \subseteq \mc{M}$. 
    That is, given a geometrically normal $G$-tiling $\mc{T}$ of $\mc{M}$, there is a $\delta$-generic Delaunay set $\mc{D} \subseteq \mc{M}$ for some $\delta > 0$ such that the simplices of $\mc{T}_{\Delta}$ are Riemannian simplices whose vertices are $\mc{D}$.
\end{proposition}

All that remains, then, is the proof of \autoref{prop:existence of generic across branches}.
That is, we must show that an indicator set on $\Gamma$ exists whose pullback to $\mc{M}$ is a $\delta$-generic Delaunay set.
The main idea is that we build up such a set inductively, by ensuring we have a Delaunay set which respects the multitude of (collared) adjacencies which a single face of a prototile admits. 
We must ensure that the subset constructed on a given prototile respects all possible adjacencies simultaneously.
That is, we require an indicator set which behaves well across branches on $\Gamma$.
To do this, we construct and adjust the indicator set inductively on neighborhoods of skeleta, starting around $0$-cells.

\begin{figure}[htb!]
    \begin{center}
    \includegraphics[width=0.45\textwidth, page=1]{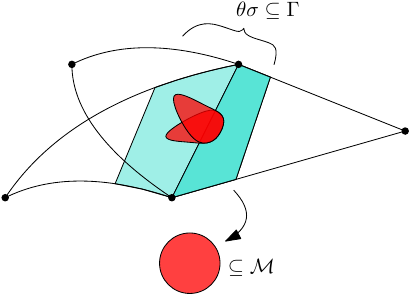}
    \hfill
    \includegraphics[width=0.45\textwidth, page=2]{Figures/consistent-across-branches.pdf}
    \end{center}
    \caption{A patch of a Delaunay set $D$ (in a subcomplex of $\Gamma$) which is $\theta$-consistent across the branch formed by the $1$-cell $\sigma$ under the coordinate charts of $\Gamma$ induced by $\mc{M}$. The corresponding ball in $\mc{M}$ serves as coordinates for both branches simultaneously.}\label{fig:consistent across branches}
\end{figure}

\begin{definition}
    Let $\pi : \mc{M} \to \Gamma$ be the natural identification for a $G$-tiling $\mc{T}$ of $\mc{M}$, and let $D \subseteq \Gamma$ be a finite indicator set.
    Let $\theta > 0$.
    We say that $D$ is $\theta$-\term{consistent across branches} if for every $k$-cell $\Sigma \subseteq \mc{M}$ of $\mc{T}$, there is a $\theta$-tubular neighborhood $\theta \Sigma$ of $\Sigma$ such that if $\pi(\Sigma) = \pi(\tilde{\Sigma})$, then $\pi^{-1}(D) \cap \tube{\theta}{\Sigma}$ and $\pi^{-1}(D) \cap \tube{\theta}{\tilde{\Sigma}}$ are identical up to translation by $G$.
    If $\theta$ is unspecified, we will simply say that $D$ is \term{consistent across branches}.
\end{definition}

Heuristically this means that even if $\Sigma$ and $\tilde{\Sigma}$ are not adjacent to the same patch of collared tiles as one another in $\mc{T}$, the point patterns around them induced by $\pi^{-1}(D)$ are identical (up to translation) because $\Sigma$ and $\tilde{\Sigma}$ are identified in $\Gamma$.
Alternatively, this means that for any ball $B_{\Gamma}$ of radius $< \theta$ centered at a branch point of $\Gamma$, there is an indicator set on a ball in $\mc{M}$ such that the natural map $B_{\Gamma} \to \mc{M}$ given by coordinates of $\Gamma$ in $\mc{M}$ will induce the patch $D \cap B_{\Gamma}$. (See \autoref{fig:consistent across branches}.)

Assuming such a consistency condition, an inductive process allows us to enumerate our indicator set in such a way that the algorithm of \autoref{thm:extended algorithm} will perturb points into a position whose pullback to $\mc{M}$ is $\delta$-generic.
In particular, by inducting up (tubular neighborhoods of) skeleta of $\Gamma$, we can ensure that the local complexity introduced by branch sets does not interfere with the perturbation strategy.
One creates a Delaunay set which is consistent across branches in the first place by a similar inductive process.

\begin{figure}[htb!]
            \begin{center}
                \includegraphics[width=0.45\textwidth, page=3]{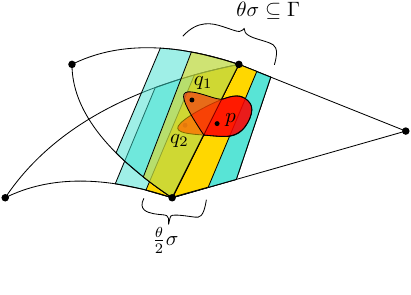}
                \hfill
                \includegraphics[width=0.45\textwidth, page=4]{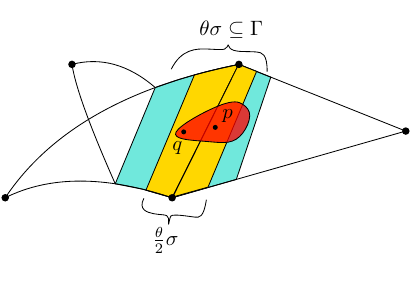}
            \end{center}
            \caption{
                Points $q_1$ and $q_2$ of a Delaunay set $D \subseteq \Gamma$ are $\theta$-consistent across the branch $\sigma$ (as in the lefthand figure).
                Because of this, the point $p$ may be perturbed into good position relative to both $q_1$ and $q_2$ simultaneously (see the righthand figure).
                The points $q_1$ and $q_2$ themselves are perturbed as if they werea single point.
                All points of $D \cap \frac{\theta}{2}\sigma$ will be perturbed into good position as this stage while those in $D \cap \theta\sigma$ more generally will be perturbed in a future inductive step.}
                \label{fig:good position across branches}
        \end{figure}

The main idea of the proof of \autoref{prop:existence of generic across branches} is that one can construct a $\delta$-generic point set, as long as we can perturb points of $\Gamma$ in ways which are consistent across branches. (See \autoref{fig:good position across branches}.)
The actual proof itself passes points back to the prototiles of the actual tiling.

\genericOnAP
    \begin{proof}
        Let $\tube{\theta}{\Gamma^{(k)}}$ denote the neighborhood of radius $\theta$ of the $k$-skeleton of $\Gamma$.
        The assumption that $\mc{T}$ is geometrically normal implies that every $k$-cell of $\mc{T}$ has a sufficiently small tubular neighborhood such that the neighborhoods of adjacent $k$-cells in $\mc{T}$ overlap only in tubular neighborhoods of lower dimensional cells.
        Hence, there is some $\theta_{k} > 0$ such that the neighborhood $\tube{\theta_{k}}\Gamma^{(k)}$ in $\Gamma$ deformation retracts to $\Gamma^{(k)}$ for each $k$.
        Fix a decreasing sequence of such radii for which $\theta_{k} \le \frac{\theta_{k-1}}{2}$.

        Inductively suppose that a sufficiently fine indicator set $D_{k-1}$ has been constructed on $\tube{\theta_{k-1}}{\Gamma^{(k-1)}}$ for which the points of $D_{k-1} \cap \tube{{\theta_{k-1}}}{\Gamma^{(k-1)}}$ are consistent across branches.
        (The base case for $k = 0$ is trivial.)
        We further assume that the points of the $\frac{\theta_{k-1}}{2}$-neighborhood of $\Gamma^{(k-1)}$ have already been perturbed into good position by the algorithm, relative to the points in the larger neighborhood $\tube{\theta_{k-1}}{\Gamma^{(k-1)}}$.
        That is, we assume the point set $D_{k-1} \cap \tube{\theta_{k-1}}{\Gamma^{(k-1)}}$ is $\delta_{k-1}$-generic for $D_{k-1} \cap \tube{\frac{\theta_{k-1}}{2}}{\Gamma^{(k-1)}}$ for some $\delta_{k-1} > 0$.
        This means that the points of $D_{k-1} \cap \tube{\frac{\theta_{k-1}}{2}}{\Gamma^{(k-1)}}$ are all $\delta_{k-1}$-protected, and induce a well-defined triangulation in a tubular neighborhood of their corresponding $(k-1)$-cells when pulled back to $\mc{M}$.
        
        We now extend $D_{k-1}$ to a set $D_{k}$ on $\tube{\theta_{k}}{\Gamma^{(k)}}$ in a way which is consistent across branches, and for which $D_k \cap \tube{\theta_{k}}{\Gamma^{(k)}}$ is $\delta_{k}$-generic for $D_{k} \cap \tube{\frac{\theta_{k}}{2}}{\Gamma^{(k)}}$ for some $\delta_{k} > 0$. 
        Let $\sigma$ be a $k$-cell of $\Gamma^{(k)}$, and let $\Sigma \subseteq \mc{T}$ denote a once-collared cell of $\mc{T}$ for which $\pi(\Sigma) = \sigma$.
        In other words, $\Sigma$ is a collared "protocell" of $\mc{T}$.
        Let $\Psi : \star(\sigma) \to \star(\Sigma)$ denote the obvious identification, and observe that $\Psi \circ \pi = \Id$ on $\star(\Sigma)$.
        (Recall that $\pi : \mc{M} \to \Gamma$ is the natural projection induced by the tiling $\mc{T}$.)

        Consider $\tube{\theta_k}{\Sigma}$, which is a tubular neighborhood of $\Sigma \subseteq \mc{M}$. Let $D_{\boundary \Sigma} = \Psi(D_{k-1} \cap \tube{\theta_{k-1}}{\boundary \sigma})$. 
        This set of points lies in the $\theta_{k-1}$ neighborhood of the boundary of $\Sigma$.
        Because $D_{k-1}$ is consistent across branches, observe that $D_{\boundary \Sigma} = \pi^{-1}(D_{k-1}) \cap \tube{\theta_{k-1}}{(\boundary \Sigma)}$.
        That is, $D_{\boundary \Sigma}$ may be identified with the restriction of $D_{k-1}$ to every coordinate sheet around $\boundary \sigma$ simultaneously.
                
        Extend $D_{\boundary \Sigma}$ arbitrarily to a sample set $D_{\tube{\theta_{k}}{\Sigma}}$ on the neighborhood $\tube{\theta_k}{\Sigma}$.
        Note that $D_{\tube{\theta_{k}}{\Sigma}}$ may be made to have arbitrarily fine mesh.
        Hence we can choose fine enough mesh to employ the tools of \autoref{thm:extended algorithm}.
        More specifically, we apply the algorithm to balls of radius $\frac{\theta_{k}}{2}$ centered at points of $D_{\tube{\theta_{k}}{\Sigma}} \cap \tube{\frac{\theta_{k}}{2}}{\Sigma}$, each in turn.
        After this perturbation process, every point of $D_{\tube{\theta_{k}}{\Sigma}} \cap \tube{\frac{\theta_{k}}{2}}{\Sigma}$ will be in good position.
        We let $D_{\tube{\theta_{k}}{\Sigma}}'$ denote the points of $D_{\tube{\theta_{k}}{\Sigma}}$ after this perturbation process is complete.

        By the choice of $\theta_{k}$ we see that $D_{\tube{\theta_{k}}{\Sigma}}'$ does not depend on any tiles outside of the immediate collar of $\Sigma$, and does not interact with $\theta_{k}$-neighborhoods of adjacent $k$-cells except in shared points of $\theta_{k-1}\boundary\Sigma$.
        In particular, the extension $D_{\tube{\theta_{k}}{\Sigma}}'$ does not depend on, or interact with, extensions to $D_{k-1}$ made to other $k$-cells.
        Therefore, $D_{\theta_{k}\sigma} := \pi(D_{\tube{\theta_{k}}{\Sigma}}')$ extends $D_{k-1}$ in a manner consistent across branches.
        Applying this process to every $k$-cell, we have an indicator set $D_{k} = \bigcup_{\sigma \subseteq \Gamma^{(k)}} D_{\theta_{k}\sigma}$ on $\theta_{k} \Gamma^{(k)}$ which is $\theta_{k}$-consistent across branches of $\Gamma^{(k)}$.
        By construction, by our application of \autoref{thm:extended algorithm}, the pullback $\pi^{-1}(D_{k})$ is $\delta_{k}$-generic for $\pi^{-1}(D_{k} \cap \tube{\frac{\theta_{k}}{2}}{\Gamma^{(k)}})$ for some $\delta_{k} > 0$.
        Because $D_{k}$ is finite and $\Gamma$ is compact, the pullback is a Delaunay set of $\mc{M}$.
        When $k = n = \dim{\mc{M}}$, we have the desired result.
    \end{proof}

    Note that in fact we have proved something stronger; not only is $\mc{D} \subseteq \mc{M}$ a $\delta$-generic Delaunay set, but $D \subseteq \Gamma$ is consistent across branches of $\Gamma$ as well.

\section{Pattern-Equivariant Cohomology}\label{sec:pattern equivariant}
    Fix a simple $G$-tiling space $\ms{T}$ of a Riemannian homogeneous space $\mc{M}$.
Our first topological result is nearly an immediate consequence of the definition, and in practice really only requires that $\ms{T}$ be geometrically normal.
\begin{theorem}[{See \cite[\S4]{sadun2003tilinginverselimits}}]\label{thm:direct limit of cohomology}
    Let $\ms{T}$ be a simple $G$-tiling space over a based Riemannian homogeneous space $\mc{M}$. Then $\check{H}^*(\ms{T}) \isom \dirlim \check{H}^*(\Gamma_n)$, where $\Gamma_n$ is the $n$th G\"{a}hler complex.
    \begin{proof}
        Recall that $\ms{T} = \invlim \Gamma_n$ by \autoref{thm:tiling spaces are inverse limits}.
        Then $\check{H}^*(\invlim \Gamma_n) \isom \dirlim \check{H}^*(\Gamma_n)$ as each $\Gamma_n$ is a finite CW complex (see \cite[Theorem 3.5]{sadun2008topology}).
    \end{proof}
\end{theorem}

In the rest of this section, we use $\Omega^k(\mc{M})$ to denote the smooth $k$-forms on a manifold $\mc{M}$. This is not to be confused with the use of $\Omega_{\mc{T}}$ to denote the hull of a tiling $\mc{T}$, and the context should always make our use of the overloaded symbol "$\Omega$" clear.
As the primary application of our definition of simple tilings, we generalize an established theorem about the topological information of tiling spaces.
Recall that a $k$-form $\omega \in \Omega^k(\mc{M})$ is said to be \term{PE to radius $r$} if, for $g \in G$ with $g.y= x$, whenever $[B_r(x)] = g.[B_r(y)]$ then $\omega_x = L_g^* \omega_y$.
We can summarize this by saying that $[B_r(x)] \sim [B_r(y)]$ implies $\omega_x \sim \omega_y$.
If $\omega$ is PE to radius $r$ for some $r > 0$, we simply say that $\omega$ is a \term{(strongly) PE $k$-form} and we denote the collection of these by $\Omega^{k}_{PE}(\mc{T})$.
The exterior derivative of a PE $k$-form will be a PE $(k+1)$-form, yielding a cochain complex
\[
    0 \to \Omega^0_{PE}(\mc{T}) \to \Omega^1_{PE}(\mc{T}) \to \cdots \to \Omega^n_{PE}(\mc{T}) \to 0
\]
and the \term{PE cohomology} of $\mc{T}$ is $H^*_{PE}(\mc{T})$, the cohomology of this complex.

One may also make sense of integer-valued PE cohomology of $\mc{T}$, though one must discuss PE cochains in place of $k$-forms.
The general ideas are nearly identical in this situation, though with some computational distinctions. 
We direct the reader to \cite{sadun2007pattern} for details on the integral version of these ideas.

The proof of \autoref{thm:PE cohomology}, though originally due to Kellendonk and Putnam in \cite{kellendonk2006ruelle} for $G = \mc{M} = \RR^n$, may be copied from \cite{sadun2007pattern} without much change, and so we refer the reader to this paper for an in-depth study. 
In the appendix to this latter paper, the following key fact is proved for a simple (polytopal) tiling of $\RR^n$:
\begin{theorem}[{\cite[Thm. 3]{sadun2007pattern}}]
    If $X$ is a branched manifold obtained by gluing (polytopal) tiles along their common boundaries, then the de Rham cohomology of $X$ is naturally isomorphic to the \u{C}ech cohomology of $X$ with real coefficients.
\end{theorem}
The proof of this requires that $X$ satisfies a Poincar\'{e} lemma. 
For tiles which are Euclidean polytopes, such a theorem is practically immediate.
Since the Poincar\'{e} lemma to be proved is a completely local result, this generalizes to Riemannian simplices.
In particular the above theorem applies to simple $G$-tilings of $\mc{M}$ with $X = \Gamma_{m}$ the $m$-G\"{a}hler complex.

\PECohomology
\begin{proof}
The proof is practically a verbatim copy of the proof in \cite{sadun2007pattern}. 
Here we provide the proof for $\RR$-coefficients.
Define the cellular map $\pi_n : \mc{M} \to \Gamma_m$ to be the obvious projection to the $m$-G\"{a}hler complex induced by the tiling $\mc{T}$.
If $[B_r(x)]_m \sim [B_r(y)]_m$ for sufficiently large $r > 0$, then clearly $\pi_m(x) = \pi_m(y)$. 
Similarly, if $\pi_m(x) = \pi_m(y)$, then there is a radius $r$ large enough so that $[B_r(x)]_m \sim [B_r(y)]_m$.

Since $\mc{T}$ is simple, $\Gamma_m$ is a sufficiently nice branched manifold so as to admit a natural notion of smooth differential $k$-forms $\Omega^k(\Gamma_m)$. It follows from the above ideas that $\omega$ is a PE $k$-form if and only if $\omega = \pi_m^*(\tilde{\omega}_{m})$ for some $m$ and some $\tilde{\omega}_{m} \in \Omega^k(\Gamma_m)$.
If $\omega$ is closed, then so is $\tilde{\omega}_{m}$, and we induce a map $\phi_m : H^*_{PE}(\mc{T}) \to H^*_{dR}(\Gamma_m;\RR)$ by $[\omega]_{PE} \mapsto [\tilde{\omega}_{m}]_{dR}$.
These fit together to form $\phi: H^*_{PE}(\mc{T}) \to \check{H}^*(\Hull{\mc{T}})$ which is well-defined, and an isomorphism.
The proof of this in \cite{sadun2007pattern} carries over as well.

To see that $\phi$ is well-defined, we must show the definition is independent of $m$ and of the representative $\omega$.
Indeed, for any $m' > m$, the inverse limit structure $\Hull{\mc{T}} = \invlim \Gamma_{m}$ yields commuting maps as in the following diagram:
\begin{diagram}
\& \Hull{\mc{T}} 
    \arrow[dl,"\rho_{m'}"'] 
    \drar{\rho_{m}} 
    \& \\
\Gamma_{m'} 
    \arrow[rr,"{\rho_{m',m}}"']
    \& \& \Gamma_{m}
\end{diagram}
Here $\rho_{m',m}$ is the forgetful map.
Hence $\tilde{\omega}_{m'} := \rho_{m',m}^*\tilde{\omega}_{m}$ and $\rho_{m'}^*\tilde{\omega}_{m'} = \rho_{m}^*\tilde{\omega}_{m}$. 
In particular then, $\tilde{\omega}_{m'}$ defines the same element of the direct limit of \v{C}ech cohomologies as $\tilde{\omega}_{m}$.

Observe that $\phi$ is injective: if $\phi([\omega]) = 0$, then there is some finite $m$ such that $[\tilde{\omega}_{m}] = 0$ as a cohomology class of $\check{H}^*(\Gamma_m)$.
Therefore $[\tilde{\omega}_{m}]_{dR} = 0$ as a cohomology class in $H^*_{dR}(\Gamma_m)$ as well.
Hence $\tilde{\omega}_{m} = d\nu$ for some $\nu \in \Omega^*(\Gamma_m)$ so that $\omega$ is exact as well as a class of $H^*_{PE}(\mc{T})$.
That is,
\[
    [\omega] = [\rho_m^*\tilde{\omega}_{m}] = [d \rho_m^*\nu] = 0.
\]

Let $\gamma \in \check{H}^*(\Omega_{\mc{T}}) = \check{H}^{*}(\invlim \Gamma_{m})$. By \autoref{thm:direct limit of cohomology} we have $\check{H}^{*}(\invlim \Gamma_{m}) \isom \dirlim \check{H}^{*}(\Gamma_{m})$ so that $\gamma$ is represented by $\gamma_m \in \check{H}^*(\Gamma_{m})$ for some finite $m$. 
Therefore, there is a closed form $\omega_m \in \Omega^*(\Gamma_m)$ representing $\gamma_m$.
Recalling that $\pi_m : \mc{T} \to \Gamma_m$ is the natural projection, we see that $\pi_{m}^* \omega_m$ is a closed PE form on $\mc{M}$, and that $\phi([\pi_m^*\omega_m]) = \gamma$.

As the G\"{a}hler complexes for a simple tiling are cellular (as in the Euclidean case), the proof for integral coefficients follows the same procedure, replacing the PE differential forms with their cellular counterpart.
\end{proof}
    
\appendix

\section{Nilpotent Lie Groups}\label{sec:nilpotent lie groups}
    This section provides background necessary to understand the progress towards \autoref{con:nilpotent tiling spaces are fiber bundles} in \autoref{sec:fiber bundles} which follows.
In particular we present background on nilpotent Lie groups, paying special attention to connected, simply connected, rational, nilpotent Lie groups with a left-invariant Riemannian metric.
The reader familiar with such groups is therefore welcome to skip this section, except perhaps as reference for future notation.

\subsection{Notation and Foundations}\label{ssec:nilpotent definitions}
    First, a \term{(real) Lie group} is a (real) smooth manifold with a compatible group structure, by which we mean that the actions of multiplication and inversion are smooth diffeomorphisms of the manifold. The tangent spaces of a Lie group $G$ can therefore all be identified by left-translation with the tangent space $\mf{g} = T_eG$ at the identity, which we call the \term{Lie algebra} of $G$. 
        The Lie algebra $\mf{g}$ associated to $G$ is a vector space together with an alternating bilinear map $[\cdot,\cdot] : \mf{g} \times \mf{g} \to \mf{g}$ called the \term{Lie bracket}. 

    A key relationship between a Lie group and its corresponding Lie algebra is given by the \term{exponential map}.
    That is, for any $x \in G$ there is a map $\exp_x : \mf{g} \to G$ for which $\exp_x(0) = x$.
        The following proposition is classical and can be found in many reasonable textbooks on the theory of Lie groups (for example, see \cite[Sec. 1.2]{corwin1990representations}).

    \begin{proposition}\label{prop:exp a diffeomorphism}
        If $G$ is a connected, simply connected nilpotent lie group, then $\exp_e : \mf{g} \to G$ is a diffeomorphism. 
    \end{proposition}

    A nilpotent group can be thought of as one which is, in some sense, "almost abelian". 
    Let $G$ be a group. 
    Recall that for $x,y \in G$, we define the commutator $[x,y] := xyx^{-1}y^{-1}$, and for subgroups $H,K$ of $G$ we define $[H,K] := \left\langle [h,k] \sothat h \in H, k \in K \right\rangle$. 
    We define $G_0 := G$, and inductively for $n > 0$, we set $G_n := [G_{n-1},G]$. 
    These subgroups fit together into a \term{lower central series}:
    \[
        G = G_0 \unrhd G_1 \unrhd G_2 \unrhd \cdots \unrhd G_n \unrhd \cdots
    \]
    In the event that $G_s = \set{e}$ for some $s \in \NN$, we say that the series \term{terminates}, and that the group $G$ is \term{nilpotent}. 
    The smallest $s$ for which $G_s = \set{e}$ is called the \term{nilpotency step} or \term{nilpotency class} of $G$. 
    Observe that if $s = 1$, the group $G$ is abelian.

    In our context, we will say that a connected, simply connected Lie group $G$ is a \term{nilpotent Lie group} if $G$ is nilpotent as a group. 
    However, we can also relate the notion of nilpotency to the Lie algebra $\mf{g}$. 
    In fact, we define $L_0 := \mf{g}$ and inductively define subspaces $L_n := [L_{n-1},L]$ where the bracket here is the Lie bracket of $\mf{g}$. 
    These again fit into a lower central series
    \[
        \mf{g} = L_0 \supseteq L_1 \supseteq L_2 \supseteq \cdots \supseteq L_n \supseteq \cdots
    \]
    In the event that $L_s = \set{0}$ for some $s \in \NN$, we say that the series \term{terminates}, and that the Lie algebra $\mf{g}$ is \term{nilpotent}. 
    The smallest $s$ for which $L_s = \set{0}$ is called the \term{nilpotency step} or \term{nilpotency class} of $\mf{g}$. 
    Observe that if $s = 1$, the Lie bracket $[\cdot,\cdot]$ is trivial, and $\mf{g} = \RR^n$ as an abelian Lie algebra.

    Often a nilpotent Lie group is defined to be one whose Lie algebra is nilpotent.
    In our context, the group (and not the algebra) structure will be our main focus, and so we take the nilpotence of $G$ for granted, stating the following as a proposition.

    \begin{proposition}
        Let $G$ be a connected, simply connected Lie group. Then $\mf{g}$ is a nilpotent Lie algebra if and only if $G$ is a nilpotent Lie group. 
    \end{proposition}
    
    The proofs of the above facts rely on what is commonly known as the \term{BCH formula}. 
    This is a way of translating between the Lie bracket of $\mf{g}$ and the group structure of $G$. 

    \begin{theorem}[BCH Formula]
        Let $X,Y \in \mf{g}$. Then $\exp(X)\cdot \exp(Y) = \exp(Z)$ where
        \[
            Z = X + Y + \frac{1}{2}[X,Y] + \frac{1}{12}\l([X,[X,Y]] + [Y,[Y,X]]\r) + \cdots 
        \]
        where "higher order" terms are successive nestings of commutators of $X$ and $Y$.
    \end{theorem}
    Note that in a nilpotent Lie group, successive commutators eventually vanish. 
    Therefore the BCH formula always terminates as a finite sum in this context.
    We note that the BCH formula has been fruitful in the study of nilpotent tilings already---for example, see \cite{kaiser2022complexity}.
    
    From now on, unless otherwise specified, all Lie groups will be assumed to be connected, simply connected, and nilpotent. 
    We will use $\GG$ (rather than $G$) to denote such a group, and though we may frequently remind the reader that $\GG$ is assumed to be nilpotent, we will not typically mention that $\GG$ is connected and simply connected.
    
    \subsection{The Heisenberg Group}\label{ssec:heisenberg}

    The first example of a non-abelian nilpotent Lie group is the $3$-dimensional \term{Heisenberg group} $\HH(\RR) = \HH^3 = \HH$. 
    This is commonly viewed as the upper triangular matrix Lie group with $1$'s along the diagonal:
    \[
        \HH = \set{\begin{pmatrix} 1 & x & z \\ 0 & 1 & y \\ 0 & 0 & 1 \end{pmatrix} \sothat x,y,z \in \RR}.
    \]
    Alternatively, we view this group as $\RR^3$ with the non-abelian group law
    \[
        (x,y,z) \star (x',y',z') = (x + x', y + y', z + xy' + z').
    \]
    This group is $2$-step nilpotent whose center $Z(\HH)$ and commutator subgroup $[\HH,\HH]$ are both the set $\set{(0,0,z) \sothat z \in \RR}$.
    It is easily seen that $\HH(\ZZ) := \set{(x,y,z) \sothat x,y,z \in \ZZ}$ and $\HH(\QQ) := \set{(x,y,z) \sothat x,y,z \in \QQ}$ are subgroups of $\HH$. 
    They are (respectively) the \term{integer} and \term{rational} Heisenberg groups. 
    Observe that in these definitions, we readily apply exponential coordinates.
    
    The Lie algebra $\mf{h}$ of $\HH$ is generated by basis $\set{X,Y,Z}$ with $\exp_e(X) = x = (1,0,0)$, and $\exp_e(Y) = y = (0,1,0)$, and $\exp_e(Z) = z = (0,0,1)$. The Lie bracket satisfies $[X,Y] = Z$, with all other brackets being trivial.
    
\subsection{Rational Lie Groups}

    In our progress towards \autoref{con:nilpotent tiling spaces are fiber bundles}, we will consider a class of Lie groups known as \term{rational} Lie groups. 

    \begin{definition}
        Let $\GG$ be a real nilpotent Lie group. We say that $\GG$ is \term{rational} if there is a basis $\set{X_1,\dots,X_n}$ of the Lie algebra $\mf{g}$ whose structure coefficients are in $\QQ$. That is, the bracket relations satisfy $[X_i,X_j] = \sum_k c_{i,j,k}X_k$ where $c_{i,j,k} \in \QQ$.
    \end{definition}

    A celebrated theorem by Malcev explains that such rational Lie groups are precisely those admitting lattices:
    
    \begin{theorem}[Malcev's Correspondence, {\cite{mal1949class}}]
        A connected, simply connected, nilpotent Lie group $\GG$ is rational if and only if $\GG$ has a lattice.
    \end{theorem}

    Here, a \term{(cocompact) lattice} of $\GG$ is a discrete subgroup $\Lambda$ for which $\GG/\Lambda$ is compact.
    The quotient of a connected, simply connected, nilpotent Lie group $\GG$ by a cocompact lattice $\Lambda$ is called a \term{nilmanifold}. 
    For example, it can be shown that the integer points $\HH(\ZZ)$ of the Heisenberg group form a lattice in $\HH$. 
    The nilmanifold $\HH/\HH(\ZZ)$ is called the \term{Heisenberg manifold}. 

    Given a rational nilpotent Lie group $\GG$, we call a lattice \term{integral} or an \term{integer lattice} when the lattice has integral exponential coordinates. 
    Such a lattice always exists for a rational nilpotent Lie group, and we will denote any such lattice by $\GG(\ZZ)$.
    Similarly, the set of \term{rational points} of $\GG$ will be denoted by $\GG(\QQ)$. 
    This is the set of points whose coordinates are rational under exponential coordinates. 

    A proof of the following "folklore" lemma about the abundance of lattices of arbitrarily small covolume in rational, nilpotent Lie groups may be found in \cite[Lemma 4.1]{Cornulier_2017}.
    \begin{proposition}
        Let $L$ be a cocompact discrete subring in $\mf{g}$ for $\GG$ a rational nilpotent Lie group. Then there are lattices $\Lambda_1,\Lambda_2 \subseteq \GG$ with $[\Lambda_2 : \Lambda_1]$ of finite index depending on $\dim(\mf{g})$ for which $\Lambda_1 \sse \exp(L) \sse \Lambda_2$.
    \end{proposition}

    As a consequence of this, we record the following corollary for later use.
    \begin{corollary}\label{cor:existence of lattices with small co-volume}
        Let $\mf{S} \sse \mf{g}$ be a finite set with rational coordinates such that $\mathrm{span}_{\RR}(\mf{S}) = \mf{g}$.
        Then there is a lattice $\Lambda \sse \GG$ such that $\gen{S} \subseteq \Lambda$, where $S := \exp(\mf{S}) := \set{\exp(s_i) \sothat s_i \in \mf{S}}$.
        \begin{proof}
            As the coordinates of $\mf{S}$ are all rational, an application of the BCH formula implies that there is a cocompact discrete subring $L$ in $\mf{g}$ containing $\mf{S}$ (say, the Lie subring generated by basis vectors scaled by the greatest common divisor of coordinates from $\mf{S}$).
            Thus $\exp(L)$ contains $S$.
            Taking $\Lambda = \Lambda_2$ from the above proposition, we see that $\gen{S} \sse \Lambda$ as claimed.
        \end{proof}
    \end{corollary}

\subsection{Left-Invariant Metrics} 
     With \autoref{prop:exp a diffeomorphism} in mind, we can (and in many ways, we should) view connected, simply connected, nilpotent Lie groups as a copy of $\RR^n$ with a non-abelian group law. 
     This perspective warrants some caution: the Euclidean metric on $\RR^n$ will not behave well with respect to a non-abelian group law.
     It is therefore natural to equip our group with a metric more naturally associated to this group structure.

    For an arbitrary group $G$, a \term{left-invariant distance} on $G$ is a metric $d : G \times G \to G$ such that, for any $a,q_1,q_2 \in G$, we satisfy $d(aq_1,aq_2) = d(q_1,q_2)$. 
    Associated to any Lie group is a family of left-invariant Riemannian metrics, each of which naturally induce a left-invariant distance. 
    In what follows, we will assume that our nilpotent Lie group has such a left-invariant Riemannian metric.

    Once a left-invariant Riemannian metric has been set on $G$, we get that $G$ is homogeneous: anything that we want to know about the local geometry of $G$ around a point $a \in G$ can be determined by considering $e \in G$ instead.
    For example, \autoref{thm:non-Euclidean normal tilings are simple} applies to tilings of $\GG$.

\section{Are Nilpotent Tiling Spaces Cantor Set Fiber Bundles?}\label{sec:fiber bundles}
    \subsection{Introduction to the Question}
In \cite{sadun2003tilinginverselimits}, it is shown that a large class of tiling spaces are inverse limits of branched manifolds.
In this process it is shown that such tiling spaces look locally like the product of a Cantor set together with a certain group $G$.
Motivated by the work of Anderson and Putnam on tiling dynamics, Sadun and Williams show in \cite{sadun2003tilingfiberbundles} that, for simple tilings of $\RR^n$, these local pieces actually fit together to create a fiber bundle over the torus $\TT^n$, providing some topological insight into these dynamical systems. 

The question is posed at the end of \cite{sadun2003tilinginverselimits} whether this fiber bundle construction extends to other spaces.
Specifically, it is asked whether the Cantor set neighborhoods constructed can be "stitched together" to yield a fiber bundle over a compact manifold.
The question is further refined by asking whether that compact manifold is the quotient of $G$ by a cocompact subgroup.
In \autoref{sec:fiber bundles} we demonstrate the following:

\begin{restatable*}
    {proposition}{nilpotentFiberBundlesWeak}\label{prop:rational implies fiber bundle}
    Suppose $\GG$ is a rational, connected, simply connected, nilpotent Lie group, and $\ms{T}$ is a rational, simple tiling space of $\GG$.
    Then there is a cocompact lattice $\Lambda \sse \GG$ such that $\ms{T}$ is a fiber bundle over the nilmanifold $\GG/\Lambda$.
\end{restatable*}

Ideally, one would like to upgrade this to a theorem which drops the "rational" assumption from the tiling space $\ms{T}$.
A proof without this extra assumption has yet to be found, as explored in \autoref{sec:fiber bundles}.
Even so, progress towards this proof demonstrates the utility of our definition of a simple tiling in the nilpotent setting. 
In fact, the proof of this fact in the Euclidean setting (see \cite{sadun2003tilingfiberbundles}) is relatively straightforward, and many of the ideas have a direct analogue in the nilpotent setting.
To this end, we here sketch the role of \autoref{thm:non-Euclidean normal tilings are simple} in making progress towards this conjecture. 
We use this sketch to illustrate some of the desirable features of a simple tiling.
We thank an anonymous referee for pointing out this heuristic argument.

By assuming that the tiling $\mc{T}$ is simple, we can essentially identify each tile with a well-defined indicator point (e.g. a center of mass).
Assume for simplicity that one of these points lies on the origin $e$ of $\GG$, and that there is a sufficiently large patch $P$ of these points around the origin so that every collared tile is represented in $P$.
A function $s$, called the \term{shape class} of the tiling, describes the displacement between indicator points.
Heuristically, the shape class $s$ describes the return vector of any indicator point $x$ to a point representing a tile equivalent to the one placed over the origin after translation.\footnote{The dual approach, which we eventually take, is defined on edges of tiles. In this case the shape class satisfies a cocycle requirement, induced by paths from the origin to any vertex of the tiling.}

The shape class of $s$ globally satisfies a finite set of conditions, as governed by $P$.
It appears that the density of rational points $\GG(\QQ)$ in $\GG$ should allow one to define a perturbation $s_{\QQ}$ of $s$ on $P$ which essentially rounds the displacement between indicator points to a nearby rational vector (i.e. an element of $\GG(\QQ)$).
This is a key challenge in the nilpotent setting, and a proof of this fact (namely, \autoref{con:arbitrarily small rational perturbations}) would allow one to drop the assumption that $\ms{T}$ is rational in \autoref{prop:rational implies fiber bundle}.

Under this assumption, one is allowed to adjust the indicator points of $P$ to a new set of points $P_{\QQ}$, such that displacements between pairs of points are now rational.
Because the tiling is simple, the induced Riemannian simplices remain stable under perturbation. 
That is, because every adjacency of tiles is represented in $P$, the combinatorial structure of adjacencies is preserved under this perturbation.
In the event that $\mc{T}$ is repetitive, we can extend $P_{\QQ}$ to a new tiling $\mc{T}_{\QQ}$ which has the same combinatorics as that of $P$.
Although $\mc{T}_{\QQ}$ is not necessarily MLD to $\mc{T}$, the two tilings are combinatorially equivalent.
The upshot is that, assuming one can find a rounding scheme such as \autoref{con:arbitrarily small rational perturbations} proposes, the hull of $\mc{T}_{\QQ}$ will be homeomorphic to that of $\mc{T}$.

Now, since the image $\im{s_{\QQ}} \subseteq \GG(\QQ)$ is finite, one can find a cocompact subgroup $\Gamma \le \GG$ which contains the subgroup generated by $\im{s_{\QQ}}$. That is to say, $\gen{\im{s_{\QQ}}} \subseteq \Gamma \le \GG$.
This ultimately allows us to pass the quotient map $\GG \to \GG/\Gamma$ to a map on the hull, which in turn yields the desired fiber bundle map $\Hull{\mc{T}} \simeq \Hull{\mc{T}_{\QQ}} \to \GG/\Gamma$.

Here we note a few of the issues of the above sketch, which help motivate our definition of a simple tiling.

First is the existence of such an indicator set; such a set is one which encodes identical patterns to those of the tiling $\mc{T}$ itself, but without any geometric distractions. 
More formally this means that we require a process for translating $\mc{T}$ into a point set $P_\mc{T}$, and then realizing $P_\mc{T}$ as a tiling $|P_\mc{T}|$, in such a way that $\mc{T}$ is MLD to $|P_\mc{T}|$.
One consequence of \autoref{thm:non-Euclidean normal tilings are simple} is that such a point pattern does exist, but it is taken through the route of a Delaunay triangulation rather than the method of a "center of mass" construction as above.
We suspect there is some analogue of a derived Vorono\"{i} tiling which uses such a construction, but our attempts to find one, or the relevant background for developing one, have so far failed.

Second, it is a subtle but crucial point that the combinatorial structure must remain stable under perturbation, and that one performs the perturbation in a pattern-equivariant manner.
That is, there are poorly selected Delaunay sets which naturally induce triangulations, but for which a perturbation of the points induces a triangulation with different combinatorics to that of the original, or (worse yet) which may not induce a single well-defined triangulation at all. 
Furthermore, if a tiling is not repetitive there is no guarantee that a patch $P$ containing representatives of every collared tile can be used to construct the desired shape function globally.

These issues are addressed by working in the dual setting (Delaunay triangulations) which has the downside of being more technical in nature.
Because of the work in \autoref{thm:non-Euclidean normal tilings are simple}, we are able to employ some powerful tools of \cite{boissonnat2013stability,boissonnat2018delaunay,dyer2015riemannian} which have recently developed conditions for algorithmically developing Delaunay sets of compact Riemannian manifolds. 
The Riemannian simplices developed are self-Delaunay, in the sense that they are the Delaunay triangulations induced by their own vertex set, and hence bypass the first issue. Moreover, their combinatorics are stable under sufficiently small perturbations, resolving the second point.
The branched manifold $\Gamma$ (i.e., the Anderson-Putnam complex) encodes the combinatorial structure of the tiling.
By working with $\Gamma$ rather than a patch $P$ of the tiling, we allow for the possibility of non-repetitive tilings.

\subsection{Assumptions and Notation}

Throughout this section, we fix a simple tiling space $\ms{T}$ of a connected, simply connected, rational, nilpotent Lie group $\GG$. 
The vertices are therefore assumed to be a $\delta$-generic Delaunay set for some $\delta > 0$.
We hereafter drop the "$\delta$" from the phrase "$\delta$-generic" and will simply say "generic" to mean "$\delta$-generic for some $\delta > 0$". 
In particular the symbol $\delta$ will be reserved throughout the rest of this paper as a coboundary operator.
The context should always make this clear.

The goal of this section is to illustrate how \autoref{thm:non-Euclidean normal tilings are simple} progresses us towards \autoref{con:nilpotent tiling spaces are fiber bundles}.
A proof of this would provide further evidence that other theorems on simple tilings of $\RR^n$ may be naturally generalized to the context of simple tilings of nilpotent Lie groups.
As we demonstrate, a proof of this conjecture rests on \autoref{con:arbitrarily small rational perturbations}, which essentially claims that if one has a finite system of equations in $\GG$ with a (nontrivial) solution, then one can find an arbitrarily close rational solution to that system.

Given a prototile $P$, we have seen that we can think of the vertex set as abstracting away the particular geometry of the tiles, allowing us to describe the combinatorics of the points in a small ball around any given vertex.
The star of any vertex in a tiling may be assumed to embed into the Anderson-Putnam complex $\Gamma$, and thus all of the relevant combinatorics and geometry will be legible from $\Gamma$ alone.
We introduce \term{shape functions} as a natural way to read this information.

\begin{definition}\label{def:shape function}
    Let $\pathspace{\Gamma}$ denote the set of paths along edges of $\Gamma^{(1)}$, and let $\loopspace{\Gamma}$ denote loops in $\Gamma^{(1)}$, based at some fixed vertex $*$.
    A \term{shape function} for $\ms{T}$ is a function $s : \pathspace{\Gamma} \to \GG$ which is a homomorphism with respect to concatenation, and which is trivial on nullhomotopic loops.
    That is, for $\gamma \in \loopspace{\Gamma}$ with $\gamma \simeq *$, we have $s(\gamma) = 1_{\GG}$.
    In fact, since $\Gamma$ is a CW complex, being a shape function is equivalent to requiring that $s(\boundary c)$ is trivial for all $2$-cells $c \subseteq \Gamma^{(2)}$.
\end{definition}

\begin{note}
    In the event that $\GG = \RR^n$, a shape function $s : \pathspace{\Gamma} \to \RR^n$ is often defined to be a $1$-cocycle.
    In the more general case, our definition certainly extends that of the abelian setting while avoiding the terminology of linear algebra presented by discussing (co)chains and (co)boundaries.
\end{note}

Observe that shape functions are defined fairly generally. There is no reference to the exact geometry of tilings in $\Gamma$, nor are there any obvious combinatorics from this tiling space in play.
We would like to have a shape function which encodes this relevant information for the tiling space $\ms{T}$, and the first task is to provide such a function.

Note also that a shape function induces a well-defined map on paths in a particular tiling $\mc{T} \subseteq \GG$.
We will often use the term "shape function" to refer either to a function on $\Gamma$ or the induced function on $\mc{T}$.
Since $\GG$ is simply connected, we see that the shape function is trivial on all (cellular) loops in $\mc{T}$, as any nullhomotopy in $\GG$ descends to a nullhomotopy in $\Gamma$ under the projection $\pi : \mc{T} \to \Gamma$.
This indicates that a good definition for the shape function associated to a tiling space $\ms{T}$ ought to be related to the displacement corresponding to an edge of a tiling.

\begin{definition}
    Let $\mathcal{E} \subseteq \GG$ be an edge of a prototile $P$ of $\ms{T}$, oriented from $x$ to $y$.
    The \term{displacement} of $\mathcal{E}$ is $d(\mc{E}) = x^{-1}y$. 
    Given an edge $E$ of $\Gamma$, we define the displacement of $E$ to be $d(E) = d(\mathcal{E})$ for any $\mc{E}$ such that $\pi(\mc{E}) = E$.
\end{definition}

Observe that this is well-defined, independent of choice of the choice of $\mc{E}$.
Indeed, by definition if $\pi(\mathcal{E}') = \pi(\mathcal{E})$, then $\mc{E}$ and $\mc{E}'$ are left-translates of one another.
In particular, suppose the edge $\mathcal{E}'$ is oriented from $u$ to $v$ in $\GG$. 
Then  there is a $g \in \GG$ for which $x = gu$ and $y = gv$ so that
\[
    d(\mathcal{E}) = x^{-1}y = (gu)^{-1}(gv) = u^{-1}v = d(\mathcal{E}').
\]

\begin{definition}
    Let $s_{\ms{T}} : \pathspace{\Gamma} \to \GG$ be defined on (oriented) edges by $s_{\ms{T}}(E) := d(E)$, and by concatenation otherwise. 
    We call $s_{\ms{T}}$ the \term{shape function associated to $\ms{T}$}.
    We might also discuss the shape function associated to a particular tiling, or to a particular prototile of a tiling under similar definitions.
\end{definition}

It is appropriate to justify such an evocative name:

\begin{proposition}
    The shape function associated to $\ms{T}$ is a shape function.
    \begin{proof}
        It is enough to check that $s$ is trivial on boundaries of $2$-cells.
        Without loss of generality, say by appropriate labeling of our tiles in the first place (or by a sufficiently fine subdivision process), we can assume that the projection $\pi : \mc{T} \to \Gamma$ is an embedding on any tile $T \sse \mc{T}$. 
        In particular, the boundary of any single $2$-cell in $\Gamma$ corresponds exactly to the boundary of a $2$-cell of some collared prototile.
        This means that if $W \subseteq \Gamma^{(2)}$ is a $2$-cell, then there is a $2$-cell $\mathcal{W} \subseteq \GG$ in some prototile $P$ for which $\pi(\mathcal{W}) = W$ is an embedding. 

        Suppose that $\boundary W = E_1 E_2 \cdots E_k$. Then $\boundary \mathcal{W} = \mathcal{E}_1 \mathcal{E}_2 \cdots \mathcal{E}_k$ where $\pi(\mathcal{E}_j) = E_j$ for each $1 \le j \le k$. 
        Let the edge $\mathcal{E}_{j}$ be directed from the vertex $v_{j-1}$ to the vertex $v_{j}$. 
        Then
        \begin{align*}
            s_{\ms{T}}(\boundary W) 
                & = s_{\ms{T}}(E_1) \cdot s_{\ms{T}}(E_2) \cdots s_{\ms{T}}(E_k) \\
                & = d(E_1) \cdot d(E_2) \cdots d(E_k) \\
                & = (v_0^{-1}v_1) \cdot (v_1^{-1} v_2) \cdot \cdots \cdot (v_{k-1}^{-1}v_k) \\
                & = v_0^{-1}v_k \\
                & = 1_{\GG}.
        \end{align*}
        The final equality holds since $\mathcal{W}$ is a $2$-cell whose boundary is $\mathcal{E}_1 \mathcal{E}_2 \cdots \mathcal{E}_k$, so that $v_k = v_0$.
        Hence $s_{\ms{T}}$ is a shape function.
    \end{proof}
\end{proposition}

\begin{note}    
    From now on, given a shape function $s$, we will often write $s(E)$ for edges $E$ in $\Gamma$ and $s(\mathcal{E})$ for corresponding edges $\mathcal{E}$ in prototiles (or tiles). 
    That is, we freely identify edges of $\mc{T}$ with edges in $\Gamma$.
\end{note} 

Suppose now that $s$ is a shape function for $\ms{T}$.
It is possible that $s$ is not the shape function associated to $\ms{T}$. 
Under what circumstances will $s$ be the shape function associated to some tiling space $\ms{T}_{s}$? 
When this is the case, how are $\ms{T}$ and $\ms{T}_{s}$ related? 
Under what circumstances will $\ms{T}$ and $\ms{T}_{s}$ be equivalent (and under which version of equivalence)?

We show below that ${s}$ will induce a new tiling space $\ms{T}_{s}$ as long as $d(s_{\ms{T}}(E),s(E))$ is small enough for every edge $E$ in $\Gamma$. 
In particular, since the Delaunay set associated to our tilings are generic, the combinatorics of the tiles are invariant under small enough perturbations.
Thus, if $d(s_{\ms{T}}(E),s(E))$ is sufficiently small, \autoref{lem:approximate tilings are homeomorphic} demonstrates that a tiling space $\ms{T}_{s}$ exists whose associated shape function is $s$. In fact $\ms{T}_{s}$ will be homeomorphic to $\ms{T}$.

Note that we cannot make $d(s_{\ms{T}}(E),s(E))$ small enough to expect that $\ms{T}$ and $\ms{T}_{s}$ are topologically conjugate, let alone MLD. 
For example, suppose that $s_{\ms{T}}$ describes the hull of a periodic tiling of $\RR^2$ by squares of unit side length, and let $s$ be arbitrarily close to $s_{\ms{T}}$, but describing squares of irrational side length. 
Even though the hulls $\ms{T}$ and $\ms{T}_{s}$ corresponding to each of these tilings are homeomorphic tori, the action of $\ZZ^2$ on the tiling space $\ms{T}$ is periodic, but would be evidently aperiodic on $\ms{T}_{s}$.

In search of a proof of \autoref{con:nilpotent tiling spaces are fiber bundles}, we first seek to "rationalize" displacements between vertices of the tiling in a consistent way that preserves the topology of the tiling space. 
This will allow us to find an underlying cocompact lattice of $\GG$ over which $\Hull{\ms{T}}$ will be a fiber bundle.

\begin{definition}
    A simple tiling (space) is called \term{rational} if the image of the associated shape function is contained in $\GG(\QQ)$.
\end{definition}

\nilpotentFiberBundlesWeak
\begin{proof}
    Since $\mc{T}$ is rational, the image of the shape function $s$ of $\mathscr{T}$ takes rational values under the inverse exponential. That is, under exponential coordinates on $\mathbb{G}$, we have 
    \[
        S \coloneq \{s(\mc{E}) \;\vert\; \mc{E} \subseteq \GG \text{ is an edge of } \mathscr{T}\} \subseteq \mathbb{G}(\QQ).
    \]
 Because the simplices of a simple tiling are nondegenerate (see \autoref{sec:delaunay}) we know that $\mathrm{Span}_{\RR} \log(S) = \mf{g}$.
 By \autoref{cor:existence of lattices with small co-volume}, there is a cocompact lattice $\Lambda \sse \GG$ containing $\langle S \rangle$ the subgroup generated by $S$.
 
 Therefore the vertices of any tiling $\mc{T} \in \ms{T}$ are equivalent mod $\Lambda$. That is, if $p : \GG \to \GG/\Lambda$ is the quotient map, then for any vertices $v_1,v_2 \in \mc{T}$ there is a path $E_1\dots E_k$ between them such that $s(E_i)  \in \gen{S} \subseteq \Lambda$. Therefore, $v_1^{-1}v_2 = s(E_1) \dots s(E_k) \in \Lambda$ and $p(v_1^{-1}v_2) = e \in \GG/\Lambda$. 

 Let $b_{\mc{T}}$ be the basepoint of $\mc{T}$. 
 The map $p_{\ms{T}} : \ms{T} \to \GG/\Lambda$ defined by $p_{\ms{T}}(\mc{T}) = p(b_\mc{T})$ is well-defined, and given $p(g) \in \GG/\Lambda$ the fiber $p_{\ms{T}}^{-1}(p(g))$ is the collection of tilings with vertices in $g\Lambda$. 
 The map $p_{\ms{T}}$ is a fiber bundle over the nilmanifold $\GG/\Lambda$. 
 Moreover, the fiber is a totally disconnected space, and as in \cite{sadun2003tilingfiberbundles} the fiber is finite if the tiling space is periodic, and otherwise is a Cantor set.
\end{proof}

Ideally, following the lead of \cite{sadun2003tilingfiberbundles}, one would like to demonstrate the following:

\begin{restatable*}{conjecture}{rationalTilings}\label{con:simple tiling spaces are homeomorphic to rational tiling spaces}
    Every simple tiling space of a rational, connected, simply connected Lie group $\GG$ is homeomorphic to a rational one.
\end{restatable*}

Assuming \autoref{con:simple tiling spaces are homeomorphic to rational tiling spaces}, we see that \autoref{prop:rational implies fiber bundle} immediately implies the following:

\begin{conjecture}\label{con:nilpotent tiling spaces are fiber bundles}
    Let $\ms{T}$ be a simple tiling space over a rational, connected, simply connected nilpotent Lie group $\GG$. 
    Then there is a cocompact lattice $\Lambda \sse \GG$ such that $\ms{T}$ is a fiber bundle over the nilmanifold $\GG/\Lambda$.
\end{conjecture}

The remainder of this section is dedicated to illustrating progress towards a potential proof of \autoref{con:simple tiling spaces are homeomorphic to rational tiling spaces}.
In order to do so, we introduce the idea of a \term{perturbation parameter}.
Informally, this is a parameter by which one may adjust the displacement between any two vertices, and still induce a combinatorially identical tiling by the induced Delaunay triangulation.
The following ideas help make this definition precise.

\begin{definition}
    Let $\varrho > 0$ and let $s_{\ms{T}}$ be a shape function for a tiling $\ms{T}$. A \term{$\varrho$-approximate associated shape function} is a shape function $s$ such that $d(s_{\ms{T}}(E),s(E)) < \varrho$ for every edge $E$ in $\Gamma$, where $\Gamma$ is the once-collared Anderson-Putnam complex. 
\end{definition}

\begin{lemma}\label{lem:perturbation parameters exist}
    Let $\ms{T}$ be a simple tiling space over a connected simply connected nilpotent Lie group $\GG$.
    Let $\Gamma$ be the once-collared Anderson-Putnam complex of $\ms{T}$.  
    Suppose $s_{\ms{T}}$ is the shape function associated to $\ms{T}$ and let $s$ be an arbitrary shape function on $\ms{T}$.
     
    There is a $\varrho > 0$ such that if $s$ is a $\varrho$-approximate associated shape function for $\ms{T}$, then for every prototile $P$ there is a simplex $P_{s}$ and a homeomorphism $\phi_{P} : P \to P_{s}$ such that $s$ is the shape function associated to $P_{s}$. 
    Moreover, the homeomorphisms on prototiles which share a common face may be made to agree. That is, if $F = P \cap g.P'$ is a $k$-simplex common to (translates of) prototiles $P$ and $P'$ of $\ms{T}$, then $\phi_{P}|_{F}$ agrees with $\phi_{P'}|_{F}$ after appropriate translation.
    \begin{proof}
        Let $P$ be a collared prototile in $\ms{T}$ with vertices $\set{p_0,p_1,\dots,p_n}$. Without loss of generality, suppose that $p_{0} = e$ is the identity. 
        Let $E_k$ denote the corresponding edge in $\Gamma$ from $e$ to $p_k$.
        Each $0$-cell $p_k$ of $P$ can be identified with an element of $\GG$ itself, we further have $p_{k} = s_{\ms{T}}(E_{k})$.

        By \autoref{prop:normal tilings are generic and simplicial}, we may assume that the vertices of $\ms{T}$ are in generic position.
        Hence, for $s$ a sufficiently small $\varrho$-perturbation, the vertex set $\set{p_0,s(E_1),s(E_2),\dots,s(E_n)}$ is still in generic position and induces a well-defined Riemannian simplex $P_{s}$ whose vertex set is $\set{p_0,s(E_1),s(E_2),\dots,s(E_n)}$. 
        Thus there is a homeomorphism $P \to P_{s}$ induced completely by the shape function $s$, which passes through the barycentric coordinate map (see \autoref{sec:delaunay}).

        To see that these homeomorphisms may be made to agree on adjacent prototiles, let $C_{P}$ denote the $0$-cells of $\star(P)$.
        Observe that $C_{P}$ is the set of vertices adjacent to $P$ in a collared version of $P$. 
        For each vertex $v \in C_{P}$ select a path $\mathcal{E}_1\mathcal{E}_2\dots \mathcal{E}_r$ in $\star(P)$ from $v$ to $p_0$. 
        Because $s$ is a shape function, the element $v^{s} = s(\mathcal{E}_1 \mathcal{E}_2 \dots \mathcal{E}_r)$ is well-defined regardless of the path chosen. 
        This then defines a vertex set $C_{P}^{s} = \set{v^{s} \sothat v \in C_P}$ which extends the vertex set of $P_{s}$. 
        Because $\ms{T}$ is a simple tiling, the set $C_P$ is in generic position, and $C_P^{s}$ is a $\varrho_P$-perturbation of $C_P$ which is in generic position for sufficiently small $\varrho_P > 0$. By \autoref{thm:stable stars of metric perturbation}, we see both that $\Del(C_P^{s})$ is a well-defined Delaunay triangulation and that in fact $$\star(P) = \Del(C_P) \cong \Del(C_P^{s}) = \star(P_{s}).$$
        In particular, $\star(P_{s})$ is a triangulation by Riemannian simplices with the same combinatorics as the collared prototile $P$ itself.
        Because the combinatorics are the same, and the homeomorphism is induced completely by local information about vertices, these homeomorphisms piece together in the desired manner. 
        Furthermore, by minimizing $\varrho = \min_{P}\set{\varrho_P}$ over all of the (finitely many) prototiles, we obtain the global quantity $\varrho > 0$ as desired. 
    \end{proof}
\end{lemma}

\begin{definition}
    The parameter $\varrho > 0$ as given in the conclusion of \autoref{lem:perturbation parameters exist} is called a \term{perturbation parameter} for $\ms{T}$.
\end{definition}

\begin{lemma}\label{lem:approximate tilings are homeomorphic}
    Let $\varrho$ be a perturbation parameter for a simple tiling space $\ms{T}$ of a connected, simply connected, nilpotent Lie group $\GG$, and suppose that $s$ is a $\varrho$-approximate associated shape function for $\ms{T}$.
    Then there is a tiling space $\ms{T}_{s}$ whose associated shape function is $s$, and for which $\ms{T}_{s}$ is homeomorphic to $\ms{T}$.
    \begin{proof}
        Fix a tiling $\mc{T} \in \ms{T}$.
        Pick an enumeration on tiles $\set{T_m}_{m \in \NN}$ on $\mc{T}$ such that $e \in T_0$, and such that $T_m \cap \bigcup_{0 \le i \le m-1} T_i \neq \emptyset$ for $m > 0$.
        Denote by $T_{[m]}$ the patch $T_{[m]} := \bigcup_{0 \le i \le m} T_i$.
        Fix $v_0 \in T_0$.

        We will create a cellular map $\phi_{\mc{T}} : \GG \to \GG$ whose image is a tiling with prototiles in $\mc{P}_{s}$, using the cellular homeomorphisms $\set{\phi_{P} \sothat P \in \mc{P}}$.
        We do this by creating a new vertex set which, relative to any vertex, looks locally like $\varrho$-perturbation. 
        By definition the new vertex set will induce a tiling whose tiles come from $\mc{P}_{s}$.

        First, let $P_0 \in \mc{P}$ and $g_0 \in P_0$ be such that $T_0 = g_0^{-1}P_0$. That is, let $g_0 \in P_0$ be the point corresponding to the placement of $e$ in the tile $T_0$.
        Define $b_0 := \phi_{P_0}(g_0)$ to be the corresponding basepoint in the $s$-prototile $P_0^{(s)}$, and define $\phi_{[0]} := L_{b_0^{-1}} \circ \phi_{P_0} \circ L_{g_0}|_{T_0}$.
        This is the map which places the prototile $P_0^{(s)}$ over the origin $e \in \GG$ in such a way that is consistent with the original placement of $e \in P_0$ in the sense that the homeomorphism induced by the barycentric coordinate map keeps $e$ fixed.

        Inductively, assume we have constructed a patch $T_{[m-1]}^{s}$ of tiles from $\mc{P}_{s}$ such that there is a cellular homeomorphism $\phi_{[m-1]} : T_{[m-1]} \to T_{[m-1]}^{s}$ for which, if $T \sse T_{[m-1]}$ is a translate of the prototile $P$, then $\phi_{[m-1]}|_{T}$ is a translate of $\phi_{P}$.

        Now, suppose that $v \in T_{m} \setminus T_{[m-1]}$ is a vertex, and let $\mathcal{E}$ be an edge from $v$ to some vertex $v_{m-1} \in T_{[m-1]}$.
        Pick a tile $T_{i}$ ($i < m$) for which $v_{m-1}$ is a vertex, so that $T_{m}$ is a tile in the collar of $T_{i}$. 
        Let $P_{m}$ be the prototile corresponding to $T_m$.
        By \autoref{lem:perturbation parameters exist}, there is a cellular homeomorphism $\zeta : [T_{i}]_{1} \to [T_{i}^{s}]_{1}$ whose restriction to $T_m$ is a translate of $\phi|_{P_{m}}$.
        By the inductive hypothesis, this cellular homeomorphism agrees with $\phi_{[m-1]}$ on adjacent tiles as well. 
        Therefore, we can extend $\phi_{[m-1]}$ by gluing on the map $\zeta$, creating a cellular homeomorphism $\phi_{[m]}$ whose image we denote by $T_{[m]}^{s}$. 
        This completes the inductive step, and hence defines a cellular homeomorphism $\phi_{\mc{T}}$ on all of $\mc{T}$ defined by $\phi_{\mc{T}}|_{T_m} := \phi_{[m]}|_{T_m}$.

        Note that this process depended on a few choices. 
        First, we chose an enumeration on tiles, and second, we chose an edge connecting the vertex $v$ to $T_{[m-1]}$.
        These choices can be summed up by saying that we chose a path $A = \mathcal{E}_1 \mathcal{E}_2 \dots \mathcal{E}_k$ from $v_0$ to each vertex $v$ where each $\mathcal{E}_i$ is an edge of a tile, and where we then define $\phi_{\mc{T}}(v) = s(A) := s(\mathcal{E}_1)s(\mathcal{E}_2) \cdots s(\mathcal{E}_k)$.
        
        We show now that this is independent of the choice of such a path $A$. 
        Since the position of the vertices uniquely determines the adjacent tile types, this shows that the resulting tiling cellular homeomorphism $\phi_{\mc{T}}$ and its image $\ms{T}_{s}$ will be independent of these choices.
        
                Fix a vertex $v_0$ in (one of) the tile(s) containing the origin, and let $v$ be any vertex in $\mc{T}$.
        Let $A$ and $B$ be two paths along edges in $\mc{T}$ from $v_0$ to $v$.
        According to the path $A$, the vertex $v$ will be placed at $s(A)$, and according to $B$, this vertex will be placed at $s(B)$.
        Observe that the path $AB^{-1}$ is a cycle in $\GG$, since $A$ and $B$ both start at $v_0$ and end at $v$. 
        Therefore, (see the discussion around \autoref{def:shape function}) we have $s(AB^{-1}) = e$ since $s$ is a shape function on $\Gamma$.
        Hence $s(A) = s(B)$, and the placement of $v$ is independent of the choice of path, and the cellular homeomorphism $\phi_{\mc{T}} : \mc{T} \to \mc{T}_{s}$ is well-defined, as claimed.

        It is clear that if $\mc{T},\mc{T}' \in \ms{T}$ are close translates of each other (say with $\mc{T}' = x\mc{T}$, for small $x$), then so are $\mc{T}_{s}$ and $\mc{T}_{s}'$ (though the exact distance will depend on the cellular homeomorphism $\phi_{P_0}$ of the origin-containing prototile).
        Similarly, if $\mc{T}$ and $\mc{T}'$ agree on a large ball around the origin, then so will $\mc{T}_{s}$ and $\mc{T}_{s}'$, as the cellular homeomorphisms will have been constructed identically on these large balls, independent of the enumeration method selected.
        Furthermore, it is clear that if $\mc{T} \neq \mc{T}'$, then $\mc{T}_{s} \neq \mc{T}_{s}'$ as well, as images of the cellular maps will disagree wherever the original tilings disagree.

        Define $\ms{T}_{s} := \set{\mc{T}_{s} \sothat \mc{T} \in \ms{T}}$. 
        Define $\phi : \ms{T} \to \ms{T}_{s}$ by $\phi(\mc{T}) = \mc{T}_{s}$. 
        By the remarks above, the cellular maps $\phi_{\mc{T}}$ work together in such a way that $\phi$ is a homeomorphism. 
        By construction, $s$ is the shape function associated to $\ms{T}_{s}$.
    \end{proof}
\end{lemma}

We now come to the crux of the argument: to prove \autoref{con:simple tiling spaces are homeomorphic to rational tiling spaces}, we must be able to find rational shape functions which approximate an associated shape function arbitrarily well.
In particular, if this is possible, we would be able to apply \autoref{lem:approximate tilings are homeomorphic} to perturb the points of a tiling into a rational approximation of the original tiling, up to a homeomorphism of the hull.

\begin{conjecture}\label{con:arbitrarily small rational perturbations}
    Let $\ms{T}$ be a simple tiling space and let $s_{\ms{T}}$ be a shape function for $\ms{T}$. There is a shape function $s^{\QQ}$ of $\ms{T}$ such that $\im{s^{\QQ}} \sse \GG(\QQ)$, and such that $d(s(E),s^{\QQ}(E))$ can be made arbitrarily small for each edge $E$ in $\Gamma$. 
\end{conjecture}

\rationalTilings
\begin{proof}[Proof (assuming {\autoref{con:arbitrarily small rational perturbations}})]
        Let $\ms{T}$ be a simple tiling space of $\GG$ with $\varrho > 0$ a perturbation parameter of $\ms{T}$.
    By \autoref{con:arbitrarily small rational perturbations}, there is a rational $\varrho$-approximate shape function $s_{\QQ}$ for $\ms{T}$.
    By \autoref{lem:approximate tilings are homeomorphic} there is a rational tiling space $\ms{T}_{\QQ}$ whose associated shape function is $s_{\QQ}$, and for which $\ms{T}_{\QQ}$ is homeomorphic to $\ms{T}$.
\end{proof}

Unfortunately, the proof of \autoref{con:arbitrarily small rational perturbations} remains evasive.
In the Euclidean setting, such a proof follows from a process of Gaussian elimination; this process does not immediately extend to the nilpotent setting, and new techniques of proof may be required.
It is clear however that a proof of \autoref{con:arbitrarily small rational perturbations} would complete a proof of \autoref{con:simple tiling spaces are homeomorphic to rational tiling spaces}, hence of \autoref{con:nilpotent tiling spaces are fiber bundles}.

\section{Further Directions}\label{sec:further directions}
    Recently, in \cite{bjorklund2018approximate}, Bj\"{o}rklund and Hartnick have introduced (uniform) approximate lattices which are closely related to mathematical quasicrystals.
In the same way that Malcev's correspondence links uniform lattices to rational nilpotent Lie groups, Machado demonstrates in \cite{machado2018approximate} that (uniform) approximate lattices are naturally associated to nilpotent Lie groups with structure constants in $\bar{\QQ} \cap \RR$. 

Such uniform approximate lattices are defined by the same methods as cut-and-project tilings, establishing a strong connection between such objects. 
A natural question, therefore, is whether there is an extension of the above results to tilings over such Lie groups.
There are two notions of (non-uniform) approximate lattices posited in \cite{bjorklund2018approximate}, both of which use the (right) hull of the approximate lattice as a substitute for the quotient $\GG/\Lambda$. 
Approximate lattices of simply connected nilpotent groups, however, are guaranteed to be uniform, and so we form the following (somewhat vague) conjecture:

\begin{conjecture}
    \autoref{con:nilpotent tiling spaces are fiber bundles} extends to nilpotent Lie groups with structure constants in $\bar{\QQ} \cap \RR$ by replacing the nilmanifold with the hull $\Hull{\Lambda}$ of $\Lambda$ in $\GG$ for some appropriately selected, "minimal" approximate lattice $\Lambda$ with displacement in $\bar{\QQ} \cap \RR$.
\end{conjecture}

What it means for $\Lambda$ to be "minimal" is not clear; perhaps algebraically this minimality can be described by some information on $\check{H}^*(\Omega_{\Lambda};\RR)$.
If there is no such algebraic minimality, perhaps there is a geometric notion of "minimality" to which $\Lambda$ subscribes. 

We also note that the results proved in the current paper have required us to use Riemannian metrics. 
For a large class of nilpotent spaces, including the Heisenberg group $\HH$, such a metric is not intrinsic. 
If there is a notion of simple tilings which allows us to include Carnot groups under their CC metrics, it would be worth investigating.
It would be especially interesting if we could extend \autoref{thm:non-Euclidean normal tilings are simple} to this context.

It is also worth observing that the rationalization process described here is not efficient. 
In the case of \autoref{ex:cut-and-project Heisenberg}, the choice of left-invariant metric on $\HH$ will change how "efficiently" the tiling covers the Heisenberg manifold, and even with a metric which looks sufficiently Euclidean at large scales, it should be apparent that the current process is inefficient.
The subdivision of tiles, rationalization of displacement vectors, as well as the construction of a fine enough viable lattice on $\HH$ all amount to making the lattice have extremely small covolume relative to the tiles.
This is true even in Euclidean space; as mentioned in \cite{sadun2003tilingfiberbundles}, applying the current process to a Penrose tiling on $\RR^2$ produces a fiber bundle which covers the torus over 20 times more than is necessary. 

When can the current process be streamlined? Are there tilings for which the current process is as efficient as possible? Answers to these questions are generally unknown, even in the case of tilings on $\RR^n$.

In a different direction of work, analysis of Euclidean tiling dynamics have been extended to the (connected, simply connected) nilpotent setting by \cite{dymarz2018separated}. 
However, there is a lack of literature investigating the study of the topology of tiling spaces in this more general setting.
As expressed earlier, we hope that the results shown here provide further evidence that such study is not only warranted, but a natural step in the life of studying tiling spaces.
We hope that our expanded notion of simple tilings will serve as a useful tool for ongoing work.


\nocite{*}
\bibliographystyle{plain}
\bibliography{sources.bib}

\end{document}